\documentclass[a4paper,11pt,leqno]{article}
\usepackage{amsmath, amssymb, amsthm}
\usepackage{graphicx} 
\usepackage[all]{xy}

\usepackage[latin1]{inputenc}

\numberwithin{equation}{section}

\begin{document}

\def\A{{\mathbb A}} \def\SSS{{\mathfrak S}} \def\BB{{\mathcal B}}
\def\PP{{\mathcal P}} \def\AA{{\mathcal A}} \def\R{{\mathbb R}}
\def\Z{{\mathbb Z}} \def\GL{{\rm GL}} \def\RR{{\mathcal R}}
\def\E{{\mathbb E}} \def\Aff{{\rm Aff}} \def\C{{\mathbb C}}
\def\Sal{{\rm Sal}} \def\B{{\mathbb B}} \def\D{{\mathbb D}}
\def\F{{\mathbb F}} \def\I{{\mathbb I}} \def\SS{{\mathcal S}}
\def\Id{{\rm Id}} \def\N{{\mathbb N}} \def\UU{{\mathcal U}}
\def\NN{{\mathcal N}} \def\H{{\mathbb H}} \def\CC{{\mathcal C}}
\def\FF{{\mathcal F}} \def\Chain{{\rm Chain}} \def\codim{{\rm codim}}
\def\Cox{{\rm Cox}} \def\S{{\mathbb S}} \def\BSal{\overline{\rm Sal}}
\def\TSal{\widetilde{\rm Sal}} \def\Ker{{\rm Ker}} \def\End{{\rm End}}

%%%%%%%%%%%%%%%%%%%%%%%%%%%%%%%%%%%%%%%%%%%%%%%%%%%%%%%%%%%%%

\title{\bf{$K(\pi,1)$ conjecture for Artin groups}}

\author{\textsc{Luis Paris}} 

\date{\today}

\maketitle

\begin{abstract}
\noindent
The purpose of this paper is to put together a large amount of results on the $K(\pi,1)$ conjecture for Artin groups, and to make them accessible to non-experts. 
Firstly, this is a survey, containing basic definitions, the main results, examples and an historical overview of the subject. 
But, it is also a reference text on the topic that contains proofs of a large part of the results on this question.
Some proofs as well as few results are new. 
Furthermore, the text, being addressed to non-experts, is as self-contained as possible. 
\end{abstract}


\section*{Introduction}

Let $X$ be a CW-complex (or a manifold having the same homotopy type as a CW-complex), and let $G$ be a discrete group. 
We say that $X$ is an Eilenberg MacLane space for $G$ if the universal cover of $X$ is contractible and its fundamental group is $G$.
From an Eilenberg MacLane space for $G$ one can easily construct a free resolution of the group algebra $\Z G$ of $G$, thus one gets a way for calculating different (co)homologies of $G$ (see \cite{Brown1}).

\bigskip\noindent
It has been proved in the 60's that the space of configurations of $n$ points in the plane is an Eilenberg MacLane space for the braid group on $n$ strands (see \cite{FoxNeu1}), and this fact has been the starting point for the calculation of the cohomology of this group \cite{Arnol1, Fuks1, Cohen1, Segal1, Vains1}.
Starting from the observation that the space of configurations of $n$ points in the plane coincides with the complement in $\C^n$ of a well-known algebraic set, the discriminant, in the 70's and 80's the construction of this space has been extended to all Artin groups as follows. 

\bigskip\noindent
By \cite{Bourb1}, it is known that any Coxeter group acts faithfully on an open nonempty convex cone $I$ so that the union of the regular orbits is the complement in $I$ of a (possibly infinite) family of linear hyperplanes. 
More generally, by \cite{Vinbe1}, if $W$ is a reflection group in Vinberg's sense (see Section 1), then $W$ is a Coxeter group, and it acts faithfully on an open nonempty convex cone $I$ so that the union of the regular orbits is the complement in $I$ of a (possibly infinite) family $\AA$ of linear hyperplanes.
By \cite{Lek1}, the fundamental group of the space
\[
N(W) = \left( (I \times I) \setminus \left( \bigcup_{H \in \AA} (H \times H) \right) \right) / W
\]
is the Artin group $A$ associated to $W$.
The $K(\pi,1)$ conjecture, due to Arnold, Brieskorn, Pham, and Thom, says that $N(W)$ is an Eilenberg MacLane space for $A$.

\bigskip\noindent
The purpose of this paper is to put together a large amount of results on this conjecture and to make them accessible to non-experts. 
Firstly, this is a survey, containing basic definitions, the main results, examples and an historical overview of the subject. But, it is also a reference text on the topic that contains proofs of a large part of the results on this question.
Some proofs as well as few results are new. 
Furthermore, the text, being addressed to non-experts, is as self-contained as possible. 

\bigskip\noindent
The paper is organized as follows.
In Section 1 we give a precise and detailed presentation of the $K(\pi,1)$ conjecture, with basic definitions, preliminaries, and examples.
Afterwards, we give an overview of the history of this question and of the cases for which the conjecture has been proved. Section 2 contains preliminaries on algebraic topology, Coxeter groups, Vinberg's reflection groups, and Artin monoids.

\bigskip\noindent
Section 3 is dedicated to a key tool of the theory: the Salvetti complexes. 
In Subsection 3.1 we define the Salvetti complex of a (possibly infinite) arrangement $\AA$ of hyperplanes in an nonempty open convex cone $I$, and we prove that this complex has the same homotopy type as the complement of $\cup_{H \in \AA} (H \times H)$ in $I \times I$.
This construction as well as the proof of this result are new, although they have been more or less known to experts. 
In Subsection 3.2 we prove that, when $\AA$ is determined by a reflection group $W$ in Vinberg's sense, then our complex coincides with the Salvetti complex defined by Charney and Davis in \cite{ChaDav2}.
In Subsection 3.3 we determine some cellular decompositions of the Salvetti complexes that we use, in particular, to show that the fundamental group of the above defined space $N(W)$ is equal to the Artin group associated to $W$.

\bigskip\noindent
In Section 4 we reprove Deligne's theorem \cite{Delig1} which says that the $K(\pi,1)$ conjecture holds if $W$ is finite.
The proof is made in a general framework in the sense that we use that $W$ is finite only in the last paragraph of the proof. 
However, we do not know how to adapt the proof in other cases. 

\bigskip\noindent
In Section 5 we study a series of results related to the $K(\pi,1)$ conjecture and to the so-called parabolic subgroups of Artin groups.
In particular, we reprove a result by Charney and Davis \cite{ChaDav1} which says that the $K(\pi,1)$ conjecture holds for Artin groups of FC type.

%%%%%%%%%%%%%%
%%CHAPITRE 1%%
%%%%%%%%%%%%%%

\section{Basic definitions, statements, and examples}

\bigskip\noindent
Let $S$ be a finite set. 
A {\it Coxeter matrix} over $S$ is a square matrix $M=(m_{s,t})_{s,t \in S}$ indexed by the elements of $S$ and satisfying (a) $m_{s,s} = 1$ for all $s \in S$; (b) $m_{s,t} = m_{t,s} \in \{2, 3, \dots\} \cup \{ \infty\}$ for all $s,t \in S$, $s \neq t$.
A Coxeter matrix is usually represented by its {\it Coxeter graph}, $\Gamma = \Gamma(M)$. 
This is a labelled graph defined as follows. 
The set of vertices of $\Gamma$ is $S$.
Two vertices $s,t \in S$ are joined by an edge if $m_{s,t} \ge 3$, and this edge is labelled by $m_{s,t}$ if $m_{s,t} \ge 4$. 

\bigskip\noindent
Let $\Gamma$ be a Coxeter graph. 
The {\it Coxeter system} of $\Gamma$ is defined to be the pair $(W,S) = (W_\Gamma,S)$, where $S$ is the set of vertices of $\Gamma$, and $W$ is the group presented as follows.
\[
W = \left\langle S\ \left|\
\begin{array}{cl}
s^2=1 & \text{for all } s \in S\\
(st)^{m_{s,t}}=1 &\text{for all } s,t \in S,\ s \neq t,\ m_{s,t} \neq \infty
\end{array}
\right. \right\rangle\,.
\]
The group $W$ is called {\it Coxeter group} of $\Gamma$.

\bigskip\noindent
{\bf Remark.}
It is shown in \cite{Bourb1} that, for $s,t \in S$, $s \neq t$, the element $st$ is of infinite order if $m_{s,t}=\infty$, and it is of order precisely $m_{s,t}$ if $m_{s,t} \neq \infty$.
Hence, the pair $(W,S)$ entirely determines the Coxeter graph $\Gamma$.

\bigskip\noindent
If $a,b$ are two letters and $m$ is an integer greater or equal to $2$, we set $\Pi(a,b : m) = (ab)^{\frac{m}{2}}$ if $m$ is even, and $\Pi(a,b : m) = (ab)^{\frac{m-1}{2}} a$ if $m$ is odd.
Let $\Sigma= \{\sigma_s ; s \in S\}$ be an abstract set in one-to-one correspondence with $S$. 
The {\it Artin system} of $\Gamma$ is defined to be the pair $(A,\Sigma)$, where $A=A_\Gamma$ is the group presented as follows.   
\[
A=\langle \Sigma \mid \Pi(\sigma_s, \sigma_t : m_{s,t}) = \Pi(\sigma_t, \sigma_s : m_{s,t}) \text{ for all } s,t \in S, s \neq t \text{ and } m_{s,t} \neq \infty \rangle\,.
\]
The group $A_\Gamma$ is called {\it Artin group} of $\Gamma$.

\bigskip\noindent
It is easily shown that the Coxeter group of $\Gamma$ admits the following presentation.
\[
W_\Gamma = \left\langle S\ \left|\
\begin{array}{cl}
s^2=1 & \text{for all } s \in S\\
\Pi(s,t : m_{s,t}) = \Pi(t,s : m_{s,t}) &\text{for all } s,t \in S,\ s \neq t,\ m_{s,t} \neq \infty
\end{array}
\right. \right\rangle\,.
\]
Hence, the map $\Sigma \to S$, $\sigma_s \mapsto s$, induces an epimorphism $\theta : A_\Gamma \to W_\Gamma$.
The kernel of $\theta$ is called {\it colored Artin group} of $\Gamma$ and it is denoted by $CA_\Gamma$.

\bigskip\noindent
{\bf Example.}
Consider the Coxeter graph $\A_n$ drawn in Figure 1.1. 
The Coxeter group of $\A_n$ has the following presentation.
\[
\left\langle s_1, \dots, s_n \ \left|\
\begin{array}{cl}
s_i^2=1 & \text{for } 1 \le i \le n\\
(s_is_{i+1})^3=1 & \text{for } 1 \le i \le n-1\\
(s_is_j)^2=1 & \text{for } |i-j| \ge 2
\end{array} \right. \right\rangle\,.
\]
This is the symmetric group $\SSS_{n+1}$ (of permutations of $\{1, \dots, n+1\}$). 
The Artin group of $\A_n$ has the following presentation. 
\[
\left\langle \sigma_1, \dots, \sigma_n \ \left|\
\begin{array}{cl}
\sigma_i \sigma_{i+1} \sigma_i = \sigma_{i+1} \sigma_i \sigma_{i+1} & \text{for } 1 \le i \le n-1\\
\sigma_i \sigma_j = \sigma_j \sigma_i & \text{for } |i-j| \ge 2
\end{array} \right. \right\rangle\,.
\]
This is the braid group $\BB_{n+1}$ on $n+1$ strands.
The colored Artin group of $\A_n$ is the pure braid group $\PP\BB_{n+1}$.

\begin{figure}[tbh]
\bigskip
\centerline{
\setlength{\unitlength}{0.5cm}
\begin{picture}(10.25,1)
\put(0,0.8){\includegraphics[width=5.2cm]{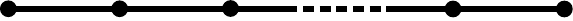}}
\put(0,0.2){\small $1$}
\put(2,0.2){\small $2$}
\put(4,0.2){\small $3$}
\put(7.5,0.2){\small $n-1$}
\put(10,0.2){\small $n$}
\end{picture}} 

\bigskip
\centerline{{\bf Figure 1.1.} The Coxeter graph $\A_n$.}
\end{figure}

\bigskip\noindent
Take a nonempty open convex cone $I$ in a finite dimensional real vector space $V$.
We define a {\it hyperplane arrangement} in $I$ to be a (possibly infinite) family $\AA$ of linear hyperplanes of $V$ satisfying (a) $H \cap I \neq \emptyset$ for all $H \in \AA$; (b) $\AA$ is {\it locally finite} in $I$, that is, for all $x \in I$, there is an open neighborhood $U_x$ of $x$ in $I$ such that the set $\{H \in \AA \mid H \cap U_x \neq \emptyset \}$ is finite.
Note that the ``classical'' definition of hyperplane arrangement imposes $I=V$ and $\AA$ finite (see \cite{OrlTer1}).

\bigskip\noindent
{\bf Example.}
Set $V=\R^3$ and $I=\{(x,y,z) \in V \mid z>0\}$.
For $k \in \Z$, we denote by $H_k$ the plane of $V$ of equation $x=kz$, and we denote by $H_k'$ the plane of equation $y=kz$.
We set $\AA=\{ H_k,\ H_k' \mid k \in \Z\}$.
This is a hyperplane arrangement in $I$.
The trace of $\AA$ on the affine plane of equation $z=1$ is represented in Figure 1.2.

\begin{figure}[tbh]
\bigskip
\centerline{
\setlength{\unitlength}{0.5cm}
\begin{picture}(9,9)
\put(1,1){\includegraphics[width=4cm]{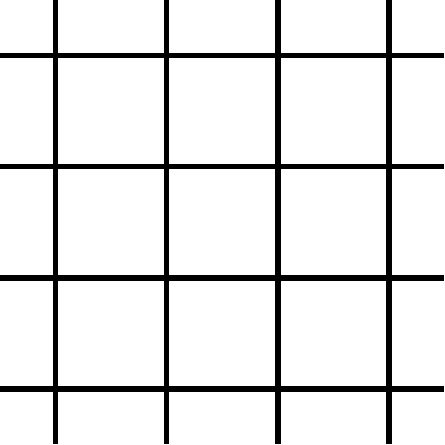}}
\put(1.4,0.2){\small $H_{-1}$}
\put(3.7,0.2){\small $H_0$}
\put(5.7,0.2){\small $H_1$}
\put(7.7,0.2){\small $H_2$}
\put(-0.5,1.8){\small $H_{-1}'$}
\put(0,3.8){\small $H_0'$}
\put(0,5.8){\small $H_1'$}
\put(0,7.8){\small $H_2'$}
\end{picture}} 

\bigskip
\centerline{{\bf Figure 1.2.} A hyperplane arrangement.}
\end{figure}

\bigskip\noindent
Let $V$ be a finite dimensional real vector space.
A {\it reflection} on $V$ is defined to be a linear transformation on $V$ of order $2$ which fixes a hyperplane.
Attention: there is no hypothesis on the orthogonality of the reflection, hence the fixed hyperplane does not necessarily determine the reflection.
Let $\bar C_0$ be a closed convex polyhedral cone in $V$ with nonempty interior, and let $C_0$ be the interior of $\bar C_0$.
A {\it wall} of $\bar C_0$ is the support of a (codimensional $1$) face of $\bar C_0$, that is, a hyperplane of $V$ generated by that face.
Let $H_1, \dots, H_n$ be the walls of $\bar C_0$.
For each $i \in \{1, \dots, n\}$ we take a reflection $s_i$ which fixes $H_i$, and we denote by $W$ the subgroup of $\GL(V)$ generated by $S=\{s_1, \dots, s_n\}$.
The pair $(W,S)$ is called a {\it Vinberg system} if $w C_0 \cap C_0 = \emptyset$ for all $w \in W \setminus \{1\}$.
In that case, the group $W$ is called {\it linear reflection group} in Vinberg's sense, $S$ is called {\it canonical generating system} for $W$, and $C_0$ is called {\it fundamental chamber} of $(W,S)$.

\bigskip\noindent
Linear reflection groups, Coxeter groups and hyperplane arrangements are linked by the following theorem.

\bigskip\noindent
{\bf Theorem 1.1}
(Vinberg \cite{Vinbe1}).
{\it Let $(W,S)$ be a Vinberg system.
We keep the above notations, and we set 
\[
\bar I = \bigcup_{w \in W} w\, \bar C_0\,.
\]
Then the following statements hold.
\begin{itemize}
\item[(1)]
$(W,S)$ is a Coxeter system. 
\item[(2)]
$\bar I$ is a convex cone with nonempty interior.
\item[(3)]
The interior $I$ of $\bar I$ is stable under the action of $W$, and $W$ acts properly discontinuously on $I$. 
\item[(4)]
Let $x \in I$ be such that $W_x=\{w \in W \mid w(x)=x\}$ is different from $\{1\}$.
Then there exists a reflection $r$ in $W$ such that $r(x)=x$.
\end{itemize}}

\bigskip\noindent
The above cone $I$ is called {\it Tits cone} of the Vinberg system $(W,S)$.

\bigskip\noindent
{\bf Remark.}
The reader must pay attention to the fact that there is a difference in Theorem~1.1.(1) between the pair $(W,S)$, viewed as a Vinberg system, and the pair $(W,S)$, viewed a Coxeter system.
Indeed, in the first case, $W$ is some specific subgroup of a linear group, while, in the second case, $W$ is just an abstract group.
Note also that any Coxeter system appears as a Vinberg system (see Theorem 2.5), but this representation is not unique in general.

\bigskip\noindent
Let $(W,S)$ be a Vinberg system.
Denote by $\RR$ the set of reflections belonging to $W$.
For $r \in \RR$ we denote by $H_r$ the fixed hyperplane of $r$, and we set $\AA=\{H_r \mid r \in \RR\}$.
Then, by Theorem~1.1, $\AA$ is a hyperplane arrangement in the Tits cone $I$.
It is called {\it Coxeter arrangement} of $(W,S)$.

\bigskip\noindent
{\bf Example.}
Consider the symmetric group $\SSS_{n+1}$ acting on the space $V=\R^{n+1}$ by permutations of the coordinates.
Let 
\[
\bar C_0 =\{x \in V \mid x_1 \le x_2 \le \cdots \le x_{n+1}\}\,.
\]
For $i,j \in \{1, \dots, n+1\}$, $i \neq j$, we denote by $H_{i,j}$ the hyperplane of equation $x_i=x_j$.
Then $\bar C_0$ is a convex polyhedral cone whose walls are $H_{1,2}, H_{2,3}, \dots, H_{n,n+1}$.
For $i \in \{1, \dots, n\}$, $s_i = (i,i+1)$ is a reflection whose fixed hyperplane is $H_{i,i+1}$.
Then $(\SSS_{n+1},\{s_1, \dots, s_n\})$ is Vinberg system. 
In this case we have
\[
\bar I = \bigcup_{w \in \SSS_{n+1}} w \bar C_0=V\,.
\]
So, $I=V$, too.
The set $\RR$ of reflections coincides with the set of transpositions, thus $\AA=\{H_{i,j} \mid 1 \le i < j \le n+1\}$ is the so-called {\it braid arrangement}.

\bigskip\noindent
{\bf Example.}
Consider the affine Euclidean plane $\E^2$. 
For $k \in \Z$, we denote by $D_k$ the affine line of equation $x=k$, and we denote by $D_k'$ the affine line of equation $y=k$ (see Figure 1.3).
We denote by $s_k$ the orthogonal affine reflection of $\E^2$ with respect to the line $D_k$, and we denote by $s_k'$ the orthogonal affine reflection with respect to $D_k'$.
We denote by $W$ the subgroup of the orthogonal affine group of $\E^2$ generated by $\{ s_k, s_k' \mid k \in \Z\}$.
We leave to the reader to determine all the elements of $W$. 
Say, however, that, among these elements, in addition to the reflections, there are U-turns, translations, and glide reflections. 
It is easily shown that $W$ is generated by $s_0, s_1, s_0', s_1'$ and admits the following presentation.
\[
W= \langle s_0, s_1, s_0', s_1' \mid s_0^2 = s_1^2 = {s_0'}^2 = {s_1'}^2=1\,,\ (s_0s_0')^2 = (s_0s_1')^2=(s_1s_0')^2 = (s_1s_1')^2 = 1\rangle\,.
\]
This is the Coxeter group of the Coxeter graph drawn in Figure 1.4.

\begin{figure}[tbh]
\bigskip
\centerline{
\setlength{\unitlength}{0.5cm}
\begin{picture}(9,9)
\put(1,1){\includegraphics[width=4cm]{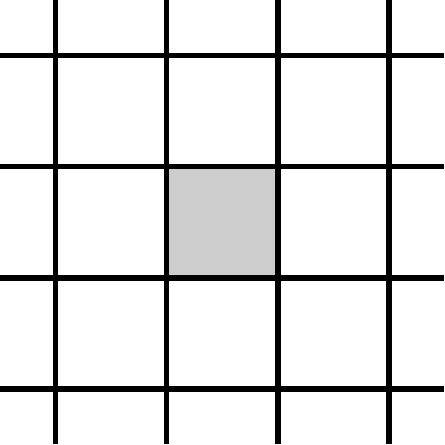}}
\put(1.4,0.2){\small $D_{-1}$}
\put(3.7,0.2){\small $D_0$}
\put(5.7,0.2){\small $D_1$}
\put(7.7,0.2){\small $D_2$}
\put(-0.5,1.8){\small $D_{-1}'$}
\put(0,3.8){\small $D_0'$}
\put(0,5.8){\small $D_1'$}
\put(0,7.8){\small $D_2'$}
\put(4.7,4.7){\small $\bar C_0'$}
\end{picture}} 

\bigskip
\centerline{{\bf Figure 1.3.} Grid lines in the affine plane.}
\end{figure}

\begin{figure}[tbh]
\bigskip
\centerline{
\setlength{\unitlength}{0.5cm}
\begin{picture}(6.4,1.5)
\put(0,0.8){\includegraphics[width=3.2cm]{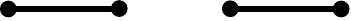}}
\put(0,0.2){\small $s_0$}
\put(2,0.2){\small $s_1$}
\put(4,0.1){\small $s_0'$}
\put(6,0.1){\small $s_1'$}
\put(0.9,1.2){\small $\infty$}
\put(4.9,1.2){\small $\infty$}
\end{picture}} 

\bigskip
\centerline{{\bf Figure 1.4.} A Coxeter graph.}
\end{figure}

\bigskip\noindent
We embed $\E^2$ in $\R^3$ via the map $(x,y) \mapsto (x,y,1)$, and we denote by $\Aff(\E^2)$ the affine group of $\E^2$.
Recall that, for all $f \in \Aff(\E^2)$, there are a unique linear transformation $f_0 \in \GL(\R^2)$ and a unique vector $u \in \R^2$ such that $f= T_u \circ f_0$, where $T_u$ denotes the translation relative to $u$.
Recall also that there is an embedding $\Aff(\E^2) \hookrightarrow \GL(\R^3)$ defined by 
\[
f \mapsto \left( 
\begin{matrix}
f_0 & u \\
0 & 1
\end{matrix}
\right)\,.
\]
Note that the elements of $\Aff(\E^2)$, embedded in $\GL(\R^3)$ via the above map, leave invariant $\E^2$ embedded into $\R^3$ as above.
So, in this way, the group $W$ can be regarded as a subgroup of $\GL(\R^3)$.
For $k \in \Z$, we denote by $H_k$ the linear plane of $\R^3$ spanned by $D_k$, and we denote by $H_k'$ the linear plane spanned by $D_k'$.
Then $s_k$ is a linear reflection whose fixed hyperplane is $H_k$, and $s_k'$ is a linear reflection whose fixed hyperplane is $H_k'$.

\bigskip\noindent
Consider the square 
\[
\bar C_0' = \{ (x,y) \in \E^2 \mid 0 \le x,y \le 1\}\,.
\]
Let $\bar C_0$ denote the cone over $\bar C_0'$.
This is a closed convex polyhedral cone whose walls are $H_0, H_1, H_0', H_1'$.
Observe that $w C_0 \cap C_0 = \emptyset$ for all $w \in W \setminus \{1\}$, thus $(W,S)$ is a Vinberg system, where $S=\{s_0,s_1, s_0', s_1'\}$.
It is easily checked that 
\[
\bar I = \bigcup_{w \in W} w\bar C_0 = \{(x,y,z) \in \R^3 \mid z >0 \} \cup \{ (0,0,0)\}\,,
\]
thus 
\[
I = \{(x,y,z) \in \R^3 \mid z >0 \}\,.
\]
On the other hand,
\[
\AA=\{ H_k, H_k' \mid k \in \Z \}\,.
\]

\bigskip\noindent
We turn now to show the link between Artin groups and Coxeter arrangements.
Besides, the $K(\pi,1)$ conjecture for Artin groups is the master peace of this link.

\bigskip\noindent
For a nonempty open convex cone $I$ in a real vector space $V$ of finite dimension $\ell$, and a hyperplane arrangement $\AA$ in $I$, we set
\[
M(\AA) = (I \times I) \setminus \left( \bigcup_{H \in \AA} H \times H \right)\,.
\]
This is a connected manifold of dimension $2\ell$.
Note that, if $I=V$, then $\AA$ is finite and 
\[
M(\AA) = (\C \otimes V) \setminus \left( \bigcup_{H \in \AA} \C \otimes H\right)\,.
\]
If $(W,S)$ is a Vinberg system and $\AA$ is the Coxeter arrangement of $(W,S)$, then we set $M(W,S) = M(\AA)$.
By Theorem 1.1, $W$ acts freely and properly discontinuously on $M(W,S)$.
Then, we set
\[
N(W,S) = M(W,S)/W\,.
\] 

\bigskip\noindent
The following result will be proved in Subsection 3.3.

\bigskip\noindent
{\bf Theorem 1.2}
(Van der Lek \cite{Lek1}).
{\it Let $(W,S)$ be a Vinberg system, and let $\Gamma$ be the Coxeter graph of the pair $(W,S)$, viewed as a Coxeter system.
Then the fundamental group of $N(W,S)$ is isomorphic to $A_\Gamma$, the fundamental group of $M(W,S)$ is isomorphic to $CA_\Gamma$, and the short exact sequence associated with the regular covering $M(W,S) \to N(W,S)$ is 
\[
\xymatrix{
1 \ar[r] & CA_\Gamma \ar[r] & A_\Gamma \ar[r]^\theta & W \ar[r] & 1}\,.
\]}

\bigskip\noindent
Recall that a space $X$ is an {\it Eilenberg MacLane space} for a discrete group $G$ if the fundamental group of $X$ is $G$ and the universal cover of $X$ is contractible.
We also say that $X$ is {\it aspherical} or that it is a {\it $K(G,1)$ space}.
Eilenberg MacLane spaces play a prominent role in cohomology of groups. 
We refer to \cite{Brown1} for more details on the subject. 

\bigskip\noindent
{\bf Conjecture 1.3}
($K(\pi,1)$ conjecture).
{\it Let $(W,S)$ be a Vinberg system, and let $\Gamma$ be the Coxeter graph of the pair $(W,S)$, viewed as a Coxeter system.
Then $N(W,S)$ is an Eilenberg MacLane space for $A_\Gamma$.}

\bigskip\noindent
Let $\AA$ be a finite hyperplane arrangement in a finite dimensional real vector space $V$.
In \cite{Salve1} Salvetti associates to $\AA$ a regular CW-complex, called {\it Salvetti complex} and denoted by $\Sal(\AA)$, and shows that $\Sal(\AA)$ has the same homotopy type as $M(\AA)$.
(The definitions of regular CW-complex and homotopy equivalence are given in Subsection 2.1.)
In Subsection 3.1 we extend the definition of $\Sal(\AA)$ to any (infinite) hyperplane arrangement $\AA$ in a nonempty open convex cone $I$, and we prove that $\Sal(\AA)$ has the same homotopy type as $M(\AA)$ (see Theorem~ 3.1).
This result is more or less known to experts, but, as far as I know, its proof does not exist anywhere in the literature.
Our proof is inspired by \cite{Paris1}.
Note that, in this paper, the complex $\Sal(\AA)$ will be defined as a simplicial complex, and, when $\AA$ is finite and $I=V$, it coincides with the barycentric subdivision of the complex originally defined by Salvetti. 

\bigskip\noindent
In Subsection 3.2, with a Coxeter graph $\Gamma$ we associate a simplicial complex $\Sal(\Gamma)$.
This complex will be naturally endowed with a free and properly discontinuous action of the Coxeter group $W$ of $\Gamma$.
Let $(W,S)$ be a Vinberg system, and let $\Gamma$ be the Coxeter graph of the pair $(W,S)$, viewed as a Coxeter system. 
We show that $\Sal(\Gamma)$ coincides with $\Sal(\AA)$, where $\AA$ is the Coxeter arrangement of $(W,S)$ (see Theorem 3.3). 
Moreover, we prove that the homotopy equivalence $\Sal(\Gamma) \to M(W,S)$ is equivariant under the action of $W$ and induces a homotopy equivalence $\Sal(\Gamma)/W \to M(W,S)/W=N(W,S)$.
In particular, this shows the following intermediate result concerning the $K(\pi,1)$ conjecture.

\bigskip\noindent
{\bf Theorem 1.4}
(Charney, Davis \cite{ChaDav1}).
{\it Let $(W,S)$ be a Vinberg system. 
Then the homotopy type of $M(W,S)$ (resp. $N(W,S)$) depends only on the Coxeter graph $\Gamma$ of the pair $(W,S)$, viewed as a Coxeter system.}

\bigskip\noindent
In their proof of Theorem 1.4, Charney and Davis \cite{ChaDav1} use another space which is homotopy equivalent to $M(W,S)$ and which depends only on the Coxeter graph $\Gamma$.
The complex $\Sal(\Gamma)$ itself is also introduced by Chaney and Davis, but in another paper \cite{ChaDav2}, and the homotopy equivalence $\Sal(\Gamma) \to M(W,S)$ is also proved in \cite{ChaDav2}.
Our proof is slightly different from the one by Charney and Davis.

\bigskip\noindent
From now on, we say that a Coxeter graph $\Gamma$ is of {\it type $K(\pi,1)$} if $\Sal(\Gamma)$ is an Eilenberg MacLane space. 
By the above, this means that $M(W,S)$ is an Eilenberg MacLane space for any representation of $(W,S)$ as a Vinberg system, where $(W,S)$ is the Coxeter system of $\Gamma$.

\bigskip\noindent
Let $\Gamma$ be a Coxeter graph, and let $(W,S)$ be the Coxeter system of $\Gamma$. 
In Subsection 3.3 we determine cellular decompositions for $\Sal(\Gamma)$ and for $\Sal(\Gamma)/W$.
The definition of $\Sal(\Gamma)$ given in Subsection~3.2 coincides with the barycentric subdivision of this cellular decomposition.
This cellular decomposition of $\Sal(\Gamma)/W$ is already defined in \cite{ChaDav2} for all Coxeter graphs, and, independently, in \cite{Salve2} when the Coxeter group $W$ is finite.
A first straightforward consequence of this description of $\Sal(\Gamma)/W$ will be that the fundamental group of $\Sal(\Gamma)/W$ (resp. $\Sal(\Gamma)$) is the Artin group $A_\Gamma$ (resp. the colored Artin group $CA_\Gamma$) (see Theorem 3.10).
This new proof of Theorem~1.2 is well-known to experts, but, as far as I know, nobody went to the bother of writing down it before. 
Note also that this cellular decomposition is a useful tool for calculating different cohomologies of $A_\Gamma$ (of course, under the condition that $\Gamma$ is of type $K(\pi,1)$) (see \cite{Calle2, Calle1, CaMoSa3, CaMoSa2, CaMoSa1,CalSal1, DePrSa1, DePrSaSt1, DeCSal2, DeCSal1, DeSaSt1, Salve2, SalStu1,Sette1,Sette2}).

\bigskip\noindent
The fact that the braid group $\BB_{n+1}$ (that is, the Artin group of $\A_n$) is of type $K(\pi,1)$ was proved in 1962 by Fox and Neuwirth \cite{FoxNeu1}.
This was the first example of an Artin group of type $K(\pi,1)$.
The $K(\pi,1)$ conjecture itself was firstly stated by Brieskorn in 1971 in \cite{Bries1}, but only for Artin groups of spherical type.
(We say that a Coxeter graph $\Gamma$ is of {\it spherical type} if the Coxeter group $W_\Gamma$ is finite.)
In the same paper, \cite{Bries1}, Brieskorn proved the conjecture for the Artin groups associated to the Coxeter graphs $\A_n$, $\B_n$, $\D_n$, $\F_4$ and $\I_2(p)$ ($p \ge 5$) (see Subsection 2.2 for the pictures of these graphs).
Immediately after, in 1972, Deligne \cite{Delig1} proved the conjecture for all spherical type Artin groups. 

\bigskip\noindent
The $K(\pi,1)$ conjecture, as is stated in the present paper, was stated for the first time in \cite{Lek1}.
According to Van der Lek, it is due to Arnold, Pham et Thom.
Besides the Artin groups of spherical type, it was previously proved in \cite{Okone1} for two families of so-called ``affine type'' Artin groups: the groups of type $\tilde \A_n$, and those of type $\tilde \C_n$ (see also \cite{ChaPei1}).

\bigskip\noindent
Let $\Gamma$ be a Coxeter graph.
For $X\subset S$, we set $M_X = (m_{s,t})_{s,t \in X}$, we denote by $\Gamma_X$ the Coxeter graph of $M_X$, and we denote by $W_X$ the subgroup of $W=W_\Gamma$ generated by $X$.
By \cite{Bourb1}, the pair $(W_X,X)$ is the Coxeter system of $\Gamma_X$.
The subgroup $W_X$ is called {\it standard parabolic subgroup} of $W$.

\bigskip\noindent
Two families of subsets of $S$ play an important role in the theory. 
The first family, denoted by $\SS^f$, consists of subsets $X \subset S$ such that $W_X$ is finite.
For $X \subset S$, we say that $\Gamma_X$ is {\it free of infinity} if $m_{s,t} \neq \infty$ for all $s,t \in X$. 
The second family, denoted by $\SS^{< \infty}$, consists of subsets $X \subset S$ such that $\Gamma_X$ is free of infinity.
Note that $\SS^f \subset \SS^{<\infty}$.

\bigskip\noindent
After \cite{Lek1}, the $K(\pi,1)$ conjecture has been proved in the following cases.
\begin{itemize}
\item[(1)]
When $m_{s,t} \ge 3$ for all $s,t \in S$, $s \neq t$ (see \cite{Hendr1}).
(Such a Coxeter graph is called of {\it large type}.)
\item[(2)]
When $|X| \le 2$ for all $X \in \SS^f$ (see \cite{ChaDav1}).
(Such a Coxeter graph is called of {\it dimension 2}.)
\item[(3)]
When $\SS^f=\SS^{<\infty}$ (see \cite{ChaDav1}).
(Such a Coxeter graph is called of {\it FC type}.)
\item[(4)]
For the ``affine type'' Artin groups of type $\tilde B_n$ (see \cite{CaMoSa1}).
\end{itemize}
Note that large type Artin groups are both, of dimension 2, and of FC type.
On the other hand, it is proved in \cite{EllSko1} that, if $\Gamma_X$ is of type $K(\pi,1)$ for all $X \in \SS^{<\infty}$, then $\Gamma$ is also of type $K(\pi,1)$ (see also \cite{GodPar1}).

\bigskip\noindent
Maybe the next advances in the subject will be due to Jon McCammond and Robert Sulway (see \cite{McCSul1}).
Indeed, they announce that they can embed any affine type Artin group into a Garside group. 
Eilenberg MacLane spaces for Garside groups are known \cite{DehLaf1,ChMeWh1}, thus such an embedding determines an Eilenberg MacLane space for the considered Artin group $A_\Gamma$.
It remains to show that this Eilenberg MacLane space has the same homotopy type as $\Sal(\Gamma)/W$. 
By the way, I thank Jon McCammond for pointing out to me this remark. 

\bigskip\noindent
In Section 4 we give a new proof of Deligne's theorem \cite{Delig1}: ``spherical type Artin groups are of type $K(\pi,1)$''.
Almost the whole proof is made for any Artin group and the hypothesis ``$\Gamma$ is of spherical type'' is used only in the last paragraph.
Nevertheless, I am not sure that this is a substantial progress toward a global proof of the $K(\pi,1)$ conjecture, as I do not know how to complete the proof for other kind of Artin groups, and Van der Lek had a similar problem (but with another complex) in \cite{Lek1} (see also \cite{Dobri1, Ozorn1}).
Our proof of Theorem 4.10 (Deligne's theorem) is inspired by the proof given in \cite{Paris1} (see also \cite{Salve3}).

\bigskip\noindent
For $X \subset S$, we set $\Sigma_X= \{\sigma_s \mid s \in X\}$, and we denote by $A_X$ the subgroup of $A=A_\Gamma$ generated by $\Sigma_X$.
Such a subgroup is called {\it standard parabolic subgroup} of $A$.
In Section 5 we use ideas from \cite{GodPar1} to prove some results that involve standard parabolic subgroups. 
In particular, we prove the following statements.
\begin{itemize}
\item[(1)]
``The pair $(A_X, \Sigma_X)$ is the Artin system of $\Gamma_X$''. 
This result is originally due to Van der Lek \cite{Lek1}.
\item[(2)]
``If $\Gamma$ is of type $K(\pi,1)$, then $\Gamma_X$ is also of type $K(\pi,1)$''.
This result, although well-known to experts, was curiously proved very recently for the first time (see \cite{GodPar1}).
\item[(3)]
``If $\Gamma_X$ is of type $K(\pi,1)$ for all $X \in \SS^{<\infty}$, then $\Gamma$ is also of type $K(\pi,1)$''.
This is the previously cited result due to Ellis and Sköldberg \cite{EllSko1}, and our proof is essentially the same as the one in \cite{EllSko1}.
Note that, thanks to Deligne's theorem \cite{Delig1}, this proves that FC type Artin groups are of type $K(\pi,1)$.
\end{itemize}

\bigskip\noindent
So, as pointed out before, many of the known results on the $K(\pi,1)$ conjecture will be proved in the present paper.
In fact, only the Artin groups of dimension 2 will not be treated, as well as some examples of Artin groups of affine type.

%%%%%%%%%%%%%%
%%CHAPITRE 2%%
%%%%%%%%%%%%%%

\section{Preliminaries}

\subsection{Preliminaries on algebraic topology}

In this subsection we present some definitions and results on algebraic topology that we will need in the sequel. 
No proof (except one) will be given, and we refer to \cite{Hatch1} for details and proofs. 

\bigskip\noindent
Let $X,Y$ be two topological spaces, and let $f,g : X \to Y$ be two continuous maps.  
We say that $f,g$ are {\it homotopic} if there exists a continuous map $H: X \times [0,1] \to Y$ such that $f(x) = H(x,0)$ and $g(x) = H(x,1)$ for all $x \in X$.
``To be homotopic'' is an equivalence relation on the set of continuous maps from $X$ to $Y$, that we denote by $\sim$.
A map $f : X \to Y$ is a {\it homotopy equivalence} if there exists a map $g: Y \to X$ such that $g \circ f \sim \Id_X$ and $f \circ g \sim \Id_Y$.
In that case we say that $X$ has the {\it same homotopy type} as $Y$.
A space $X$ is {\it contractible} if it has the same homotopy type as a point. 

\bigskip\noindent
We say that a subspace $Y$ of a topological space $X$ is a {\it deformation retract} of $X$ if there exists a continuous map $H : X \times [0,1] \to X$ such that $H(x,0)=x$ and $H(x,1) \in Y$ for all $x \in X$, and $H(y,t) =y$ for all $(y,t) \in Y \times [0,1]$.
Clearly, if $Y$ is a deformation retract of $X$, then the inclusion $Y \to X$ is a homotopy equivalence. 
The reverse is true when $X$ is a CW-complex and $Y$ is a subcomplex of $X$ (see \cite[Thm. 4.5]{Hatch1}).

\bigskip\noindent
A {\it CW-complex} is defined to be a topological space $X$ endowed with a filtration by closed subspaces,
\[
X_0 \subset X_1 \subset \cdots \subset X_n \subset X_{n+1} \subset \cdots
\]
satisfying the following properties.
\begin{itemize}
\item[(a)]
$X_0$ is a discrete set.
\item[(b)]
For all $n \in \N$, there exists a collection $\BB_n$ of $n$-dimensional closed disks, called {\it $n$-cells}, and, for each $\B \in \BB_n$, there exists a map $\varphi_\B : \partial\B \to X_{n-1}$, such that the image of each $\varphi_\B$ is a finite union of cells of $X_{n-1}$, and $X_n$ is obtained from $X_{n-1}$ gluing each $\B \in \BB_n$ to $X_{n-1}$ via the map $\varphi_\B$.
\item[(c)]
$X=\cup_{n=0}^\infty X_n$, and $X$ is endowed with the {\it weak topology} (that is, $A \subset X$ is closed if and only if $A \cap X_n$ is closed for all $n \in \N$).
\end{itemize}
We say, moreover, that $X$ is {\it regular} if, for all $n \in \N$ and all $\B \in \BB_n$, the gluing map $\varphi_\B: \partial \B \to X_{n-1}$ is a homeomorphism onto its image.  
For $n \in \N$, the subspace $X_n$ is called {\it $n$-skeleton} of $X$.

\bigskip\noindent
Let $X'$ be another CW-complex.
We denote by $\BB_n'$ the set of $n$-cells of $X'$ and, for $\B \in \BB_n'$, we denote by $\varphi_\B' : \partial \B \to X_{n-1}'$ the gluing map.
We say that $X'$ is a {\it subcomplex} of $X$ if, for all $n \in \N$, $\BB_n'$ is included in $\BB_n$, and, for all $\B \in \BB_n'$, the map $\varphi_\B' : \partial \B \to X_{n-1}'$ coincides with $\varphi_\B: \partial \B \to X_{n-1}$.

\bigskip\noindent
An (abstract) {\it simplicial complex} is defined to be a pair $\Upsilon = (S,A)$, where $S$ is a set, called {\it set of vertices}, and $A$ is a set of subsets of $S$, called {\it set of simplices}, satisfying the following properties.
(a) $\emptyset$ is not a simplex, and all the simplices are finite.
(b) All the singletons are simplices.
(c) Any nonempty subset of a simplex is a simplex.

\bigskip\noindent
Let $\Upsilon=(S,A)$ be a simplicial complex.
Take an abstract set $B=\{e_s \mid s \in S\}$ in one-to-one correspondence with $S$, and denote by $V$ the real vector space having $B$ as a basis.
For $\Delta=\{s_0,s_1, \dots, s_p\}$ in $A$, we set
\[
|\Delta| = \{t_0e_{s_0}+t_1e_{s_1}+ \cdots + t_p e_{s_p} \mid 0 \le t_0,t_1, \dots, t_p \le 1 \text{ and } \sum_{i=0}^p t_i=1 \}\,.
\]
Note that $|\Delta|$ is a (geometric) simplex of dimension $p$.
In particular, $|\Delta|$ is topologically a $p$-dimensional disk.
The {\it geometric realization} of $\Upsilon$ is defined to be the following subset of $V$.
\[
|\Upsilon| = \bigcup_{\Delta \in A} |\Delta|\,.
\]
We endow $|\Upsilon|$ with the weak topology (see \cite[Chap. III, Sec. 1]{Spani1}), so that $|\Upsilon|$ is naturally endowed with a structure of regular CW-complex. 
For $p \in \N$, if $\Delta$ is a simplex of cardinality $p+1$, then $|\Delta|$ is a cell of dimension $p$.
The geometric realization of an abstract simplicial complex is called {\it geometric simplicial complex}.

\bigskip\noindent
If $(E,\le)$ is a partially ordered set, then the nonempty finite chains of $E$ form a simplicial complex, called {\it derived complex} of $(E, \le)$ and denoted by $E'=(E,\le)'$.
This observation is of importance in the paper as our different versions of the Salvetti complex will be defined as geometric realizations of derived complexes of ordered sets. 

\bigskip\noindent
Let $X$ be regular a CW-complex.
Denote by $\BB$ the set of all cells of $X$.
If $\B$ and $\B'$ are two cells of $X$ of dimension $n$ and $m$, respectively, such that $n <m$ and $\B \subset \varphi_{\B'}(\partial \B')$, then we set $\B < \B'$.
It is easily checked that the relation $\le$ on $\BB$, defined by $\B \le \B'$ if either $\B<\B'$ or $\B=\B'$, is a partial order relation. 
The derived complex of $(\BB, \le)$ is called {\it barycentric subdivision} of $X$.
It is easily shown that $|(\BB,\le)'|$ is homeomorphic to $X$.

\bigskip\noindent
Let $X$ be a topological space, and let $\UU$ be a cover of $X$ by open subsets. 
The {\it nerve} of $\UU$, denoted by $N(\UU)$, is the simplicial complex defined as follows.
(a) The vertices of $N(\UU)$ are the elements of $\UU$.
(b) A nonempty finite set of vertices $\{U_0,U_1, \dots, U_p\}$ is a simplex in $N(\UU)$ if $U_0 \cap U_1 \cap \cdots \cap U_p \neq \emptyset$.

\bigskip\noindent 
We leave to the reader to look in the literature for the definition of a paracompact space.
However, we point out that all the spaces that we will consider are paracompact. 
The following result is one of the main tools in the paper.
Its proof can be found for instance in \cite[Sec. 4G]{Hatch1}.

\bigskip\noindent
{\bf Theorem 2.1.}
{\it 
\begin{itemize}
\item[(1)]
Let $X$ be a paracompact space, and let $\UU$ be a cover by open subspaces such that any finite nonempty intersection of elements of $\UU$ is contractible. 
Then the geometric realization $|N(\UU)|$ of the nerve of $\UU$ is homotopy equivalent to $X$.
\item[(2)]
Let $X$ be a CW-complex. 
Suppose there exists an infinite chain
\[
Y_0 \subset Y_1 \subset Y_2 \subset \cdots Y_n \subset Y_{n+1} \subset \cdots
\]
of subcomplexes of $X$ such that $Y_n$ is contractible for all $n \in \N$, and $\cup_{n \in \N} Y_n=X$.
Then $X$ is contractible, too.
\end{itemize}}

\bigskip\noindent
At some point, we will need an equivariant version of Theorem 2.1.(1), and, for this, we will need an explicit description of the homotopy equivalence $|N(\UU)| \to X$.

\bigskip\noindent
We take a connected paracompact space $X$ and a cover $\UU$ of $X$ by open subsets such that every nonempty finite intersection of elements of $\UU$ is contractible.
We denote by $\NN\UU$ the set of all finite nonempty intersections of elements of $\UU$ ordered by the inclusion.
One can show (with some effort) that $|\NN\UU'|=|N(\UU)|$, where $N(\UU)$ denotes the nerve of $\UU$.
We describe the homotopy equivalence $f: |\NN\UU'| \to X$ on the $n$-skeleton of $|\NN\UU'|$ by induction on $n$.

\bigskip\noindent
Let $U \in \NN\UU$.
Denote by $\delta(U)$ the vertex of $| \NN\UU'|$ corresponding to $U$.
Choose a point $x \in U$, and set $f(\delta(U))=x$.
This defines $f : |\NN\UU'|_0 \to X$.
Let $U_0 \subset U_1$ be a chain of length $2$ in $\NN\UU$, and let $\Delta = \Delta(U_0,U_1)$ be the $1$-simplex of $|\NN\UU'|$ corresponding to this chain.
By construction, $U_0 \subset U_1$, and, by hypothesis, $U_1$ is connected, thus there exists a path $\gamma : [0,1] \to U_1$ such that $\gamma(0) = f(\delta(U_0))$ and $\gamma(1) = f(\delta(U_1))$.
We define $f: \Delta \to U_1 \subset X$ by
\[
f((1-t)\delta(U_0) + t \delta (U_1)) = \gamma(t)\,, \quad t \in [0,1]\,.
\]
This defines the map $f : |\NN\UU' |_1 \to X$.

\bigskip\noindent
We assume that $n \ge 1$ and that the map $f: | \NN\UU' |_n \to X$ is constructed.
Furthermore, we assume that, if $U_0 \subset U_1 \subset \dots \subset U_n$ is a chain of length $n+1$ and $\Delta = \Delta(U_0,U_1, \dots, U_n)$ is the corresponding $n$-simplex of $|\NN\UU'|$, then $f(\Delta) \subset U_n$.
Let $U_0 \subset U_1 \subset \dots \subset U_{n+1}$ be a chain of length $n+2$ in $\NN\UU$, and let $\Delta = \Delta(U_0,U_1, \dots, U_{n+1})$ be the corresponding $(n+1)$-simplex in $| \NN\UU'|$.
By the above, we have $f(\partial \Delta) \subset U_{n+1}$ and, by hypothesis, $U_{n+1}$ is contractible, thus $f|_{\partial \Delta}$ extends to a continuous map $f: \Delta \to U_{n+1} \subset X$.
This defines $f : | \NN\UU' |_{n+1} \to X$.

\bigskip\noindent
The following result is probably known, but I have not found it in the literature, thus I include a proof. 
It is of importance to prove that the homotopy equivalence $\Sal(\Gamma) \to M(W,S)$ induces a homotopy equivalence $\Sal(\Gamma)/W \to M(W,S)/W= N(W,S)$ (see Corollary 3.4).

\bigskip\noindent
{\bf Proposition 2.2.}
{\it Let $X$ be a paracompact space, and let $\UU$ be a cover by open subsets such that every finite nonempty intersection of elements of $\UU$ is contractible.
Let $G$ be a group acting freely and properly discontinuously on $X$, and such that, for all $g \in G \setminus \{1\}$ and all $U \in \UU$, we have $g(U) \in \UU$ and $U \cap g(U) = \emptyset$.
Then $G$ acts freely and properly discontinuously on $|N(\UU)|$, there exists a $G$-equivariant homotopy equivalence $f : |N(\UU)| \to X$, and this homotopy equivalence induces a homotopy equivalence $\bar f : | N(\UU)|/G \to X/G$.}

\bigskip\noindent
{\bf Proof.}
By construction, the group $G$ acts on $\UU$, and sends every simplex of $|\NN\UU'|$ to a simplex, thus this action induces an action of $G$ on $|N(\UU)|=|\NN\UU'|$.
It is easily checked that the latter action is free and properly discontinuous.
On the other hand, it is easily seen that the homotopy equivalence $f: |\NN\UU'| \to X$ described above can be made to be equivariant under the actions of $G$. 
In particular, such a $f$ induces a continuous map $\bar f : |\NN\UU'|/G \to X/G$.
It remains to show that $\bar f$ is a homotopy equivalence. 

\bigskip\noindent
In order to prove that $\bar f$ is a homotopy equivalence, we will use the following results. 
These are classical and well-known.
We refer to \cite[Chap. 4]{Hatch1} for their proofs. 
Let $\varphi: X \to Y$ be a continuous map between two connected spaces having the homotopy type of CW-complexes.
Let $x_0 \in X$ be a base point for $X$, and set $y_0 = \varphi(x_0)$.
\begin{itemize}
\item[(1)]
 If $\varphi$ is a covering map, then $\varphi$ induces an isomorphism $\varphi_*: \pi_n(X,x_0) \to \pi_n(Y,y_0)$ for all $n \ge 2$.
\item[(2)]
If $\varphi$ is a regular covering map, and $G$ is its Galois group, then we have a short exact sequence 
\[
\xymatrix{
1 \ar[r] & \pi_1(X,x_0) \ar[r]^{\varphi_*} & \pi_1(Y,y_0) \ar[r] & G \ar[r] & 1}\,.
\]
\item[(3)]
The map $\varphi: X \to Y$ is a homotopy equivalence if and only if the homomorphism $\varphi_*: \pi_n(X,x_0) \to \pi_n(Y,y_0)$ is an isomorphism for all $n \ge 1$.
\end{itemize}
The group $\pi_n(X,x_0)$ is the $n$-th homotopy group of $X$. 
The reader do not need its definition to understand the proof of Proposition 2.2.
The above properties suffice.

\bigskip\noindent
We denote by $p:| \NN\UU' | \to |\NN\UU'|/G$ and by $p': X \to X/G$ the natural projections.
For $n \ge 2$ we have the following commutative diagram. 
\[
\xymatrix{
\pi_n(| \NN\UU' |) \ar[r]^{f_*} \ar[d]^{p_*} & \pi_n(X) \ar[d]^{p'_*} \\
\pi_n(| \NN\UU' |/G) \ar[r]^{\bar f_*} & \pi_n(X/G)}
\]
By (1), $p_*$ and $p'_*$ are isomorphisms, and, by (3), $f_*$ is an isomorphism, thus $\bar f_*$ is an isomorphism, too.
For $n=1$ we have the following commutative diagram, where the rows are exact sequences.
\[
\xymatrix{
1 \ar[r] & \pi_1(| \NN\UU' |) \ar[r]^{p_*} \ar[d]^{f_*} & \pi_1(| \NN\UU' |/G) \ar[r] \ar[d]^{\bar f_*} & G \ar[d]^{\Id} \ar[r] & 1\\
1 \ar[r] & \pi_1(X) \ar[r]^{p'_*} & \pi_1(X/G) \ar[r] & G \ar[r] &1}
\]
Since $f_*$ is an isomorphism, by the five lemma, $\bar f_*$ is an isomorphism, too.
By (3) we conclude that $\bar f$ is a homotopy equivalence.
\qed

\bigskip\noindent
In order to show that the fundamental group of $\Sal(\Gamma)/W$ is the Artin group $A_\Gamma$ (see Theorem~3.10), we will need the following method for computing fundamental groups of CW-complexes. 

\bigskip\noindent
Take a connected CW-complex $X$.
As in the definition, for $n \in \N$, we denote by $\BB_n$ the set of $n$-dimensional cells of $X$, and, for $\B \in \BB_n$, we denote by $\varphi_\B : \partial \B \to X_{n-1}$ the gluing map of $\B$.
Let $a \in \BB_1$ be a $1$-cell.
We set an orientation on $a$.
This means that we choose some identification of $a$ with the interval $[0,1]$.
In that way, $a$ determines a path $\tilde a : [0,1] \to X_1$ by setting $\tilde a(0) = \varphi_a(0)$, $\tilde a(1) = \varphi_a(1)$, and $\tilde a(t) = t$ for all $t \in (0,1)$.
Let $\B \in \BB_2$ be a $2$-cell.
Then $\B$ is homeomorphic to the disk $\D=\{ z \in \C \mid |z| \le 1\}$.
Without loss of generality, we can assume that $\varphi_\B(1)$ is a vertex $x_0 \in X_0$.
Then the map $\tilde \varphi_\B : [0,1] \to X_1$ defined by
\[
\tilde \varphi_\B(t) = \varphi_\B(e^{2i\pi t})
\]
is a loop based at $x_0$ homotopic in $X_1$ to a loop at $x_0$ of the form $\tilde a_1^{\varepsilon_1} \cdots \tilde a_\ell^{\varepsilon_\ell}$, with $a_1, \dots, a_\ell \in \BB_1$, and $\varepsilon_1, \dots, \varepsilon_\ell \in \{ \pm 1\}$.
Recall finally that a {\it maximal tree} of the $1$-skeleton $X_1$ is a subcomplex $T$ of $X_1$ such that $T_0=X_0$, and $T$ is simply connected. 

\bigskip\noindent
Fix a maximal tree $T$ of $X_1$ and a base-point $x_0 \in X_0$.
For all $x \in X_0$, choose a path $\gamma_x$ in $T$ from $x_0$ to $x$.
Note that $\gamma_x$ is unique up to homotopy, since $T$ is simply connected.
For a loop $\alpha : [0,1] \to X$ based at $x_0$ we denote by $[\alpha]$ the element of $\pi_1(X,x_0)$ represented by $\alpha$.
For $a \in \BB_1$ we set
\[
s_a = [\gamma_{\tilde a(1)}^{-1}\, \tilde a\,\gamma_{\tilde a(0)}]\,.
\]
Note that, if $a$ is a $1$-cell of $T$, then $s_a=1$ in $\pi_1(X,x_0)$.
On the other hand, for $\B \in \BB_2$, we take a loop of the form $\tilde a_1^{\varepsilon_1} \cdots \tilde a_\ell^{\varepsilon_\ell}$ based at $\varphi_\B(1)$ and homotopic in $X_1$ to $\tilde \varphi_\B$, and we set
\[
w(\B) = s_{a_1}^{\varepsilon_1} \cdots s_{a_\ell}^{\varepsilon_\ell}\,.
\]
Note that we have $w(\B)=1$ in $\pi_1(X,x_0)$ for all $\B \in \BB_2$.
The following result is classical in the study of CW-complexes. 

\bigskip\noindent
{\bf Theorem 2.3.}
{\it Take a connected CW-complex $X$, and keep the above notations.
Then $\pi_1(X,x_0)$ has a presentation with generators $s_a$, $a \in \BB_1$, and relations
\begin{gather*}
s_a=1 \quad \text{for all edges } a \text{ of } T\,,\\
w(\B)=1 \quad \text{for all } \B \in \BB_2\,.
\end{gather*}}

\bigskip\noindent
The proof of the following is contained in the proof of \cite[II, Thm. 7.3]{Brown1}.
It will be the key tool in the proof of Theorem 5.6.

\bigskip\noindent
{\bf Theorem 2.4.}
{\it Let $X$ be a CW-complex which is the union of two subcomplexes, $X_1$ and $X_2$, whose intersection, $Y$, is nonempty and connected.
We take a basepoint $x_0 \in Y$, and we denote by $\iota_i : \pi_1(Y,x_0) \to \pi_1(X_i,x_0)$ the homomorphism induced by the inclusion $Y \hookrightarrow X_i$, for $i=1,2$.
We assume that
\begin{itemize}
\item[(a)]
$\iota_1$ and $\iota_2$ are injective,
\item[(b)]
$X_1$, $X_2$, and $Y$ are Eilenberg MacLane spaces.
\end{itemize}
Then $X$ is also an Eilenberg MacLane space.}

\subsection{Preliminaries on Coxeter groups}

\bigskip\noindent
Let $\Gamma$ be a Coxeter graph, and let $(W,S)$ be its Coxeter system.
Take an abstract set $\{e_s \mid s \in S\}$ in one-to-one correspondence with $S$, and denote by $V$ the real vector space having $\{e_s \mid s \in S\}$ as a basis.
Define the symmetric bilinear form $B: V \times V \to \R$ by 
\[
B(e_s,e_t) = \left\{ \begin{array}{ll}
- \cos (\frac{\pi}{m_{s,t}}) & \text{if } m_{s,t} \neq \infty\\
-1 & \text{if } m_{s,t} = \infty
\end{array}\right.
\]
For $s \in S$ define $\rho_s \in \GL(V)$ by  
\[
\rho_s(x) = x - 2\, B(x,e_s)e_s\,, \quad x \in V\,.
\]
Then $\rho_s$ is a linear reflection for all $s \in S$, and the map $S \to \GL(V)$, $s \mapsto \rho_s$, induces a linear representation $\rho: W \to \GL(V)$ (see \cite{Bourb1}).
This linear representation is called {\it canonical representation} of $(W,S)$.

\bigskip\noindent
Denote by $V^*$ the dual space of $V$.
Recall that any linear map $f \in \GL(V)$ determines a linear map $f^t \in \GL(V^*)$ defined by 
\[
\langle f^t(\alpha), x \rangle = \langle \alpha, f(x) \rangle\,,
\]
for all $\alpha \in V^*$ and all $x \in V$.
The {\it dual representation} $\rho^* : W \to \GL(V^*)$ of $\rho$ is defined by 
\[
\rho^*(w) = (\rho(w)^t)^{-1}\,,
\]
for all $w \in W$.
For $s \in S$, we set $H_s = \{ \alpha \in V^* \mid \langle \alpha, e_s \rangle = 0 \}$. 
Let 
\[
\bar C_0 = \{ \alpha \in V^* \mid \langle \alpha, e_s \rangle \ge 0 \text{ for all } s \in S\}\,.
\]

\bigskip\noindent
{\bf Theorem 2.5}
(Tits \cite{Tits1}, Bourbaki \cite{Bourb1}).
{\it Let $\Gamma$ be a Coxeter graph, and let $(W,S)$ be its Coxeter system.
\begin{itemize}
\item[(1)]
The canonical representation $\rho : W \to \GL(V)$ and the dual representation $\rho^*: W ^* \to \GL(V^*)$ are faithful.
\item[(2)]
The set $\bar C_0$ is a simplicial cone whose walls are $H_s$, $s \in S$.
The transformation $\rho^*(s)$ is a linear reflection whose fixed hyperplane is $H_s$, for all $s \in S$.
Moreover, we have $\rho^*(w) C_0 \cap C_0 = \emptyset$ for all $w \in W \setminus \{1\}$.
\end{itemize}
In particular, $(\rho^*(W), \rho^*(S))$ is a Vinberg system whose associated Coxeter graph is $\Gamma$.}

\bigskip\noindent
Recall that $\Gamma$ (resp. $A_\Gamma$) is said to be of {\it spherical type} if $W_\Gamma$ is finite.
Note that, if $\Gamma_1, \dots, \Gamma_\ell$ are the connected components of $\Gamma$, then $W_\Gamma = W_{\Gamma_1} \times \cdots \times W_{\Gamma_\ell}$.
In particular, $\Gamma$ is of spherical type if and only if all its connected components are of spherical type. 

\bigskip\noindent
{\bf Theorem 2.6}
(Coxeter \cite{Coxet1, Coxet2}).
{\it
\begin{itemize}
\item[(1)]
The Coxeter graph $\Gamma$ is of spherical type if and only if the bilinear form $B : V \times V \to \R$ is positive definite.
\item[(2)]
The spherical type connected Coxeter graphs are precisely those listed in Figure 2.1.
\end{itemize}}

\begin{figure}[tbh]
\bigskip
\centerline{
\setlength{\unitlength}{0.5cm}
\begin{picture}(24,24.9)
\put(2,0){\includegraphics[width=11cm]{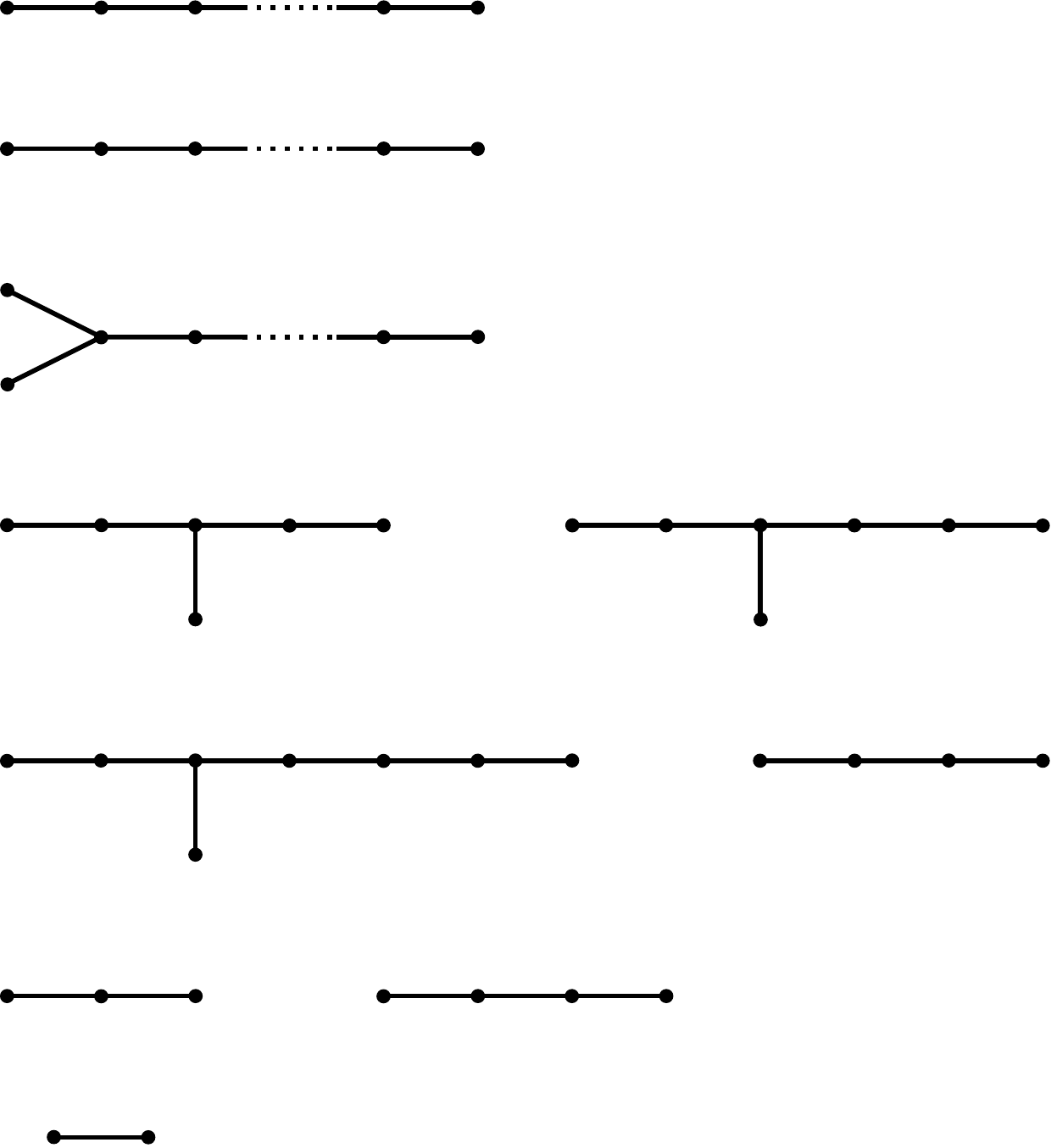}}
\put(0.5,23.7){\small $\A_n$}
\put(13,23.8){\small $n\ge 1$}
\put(0.5,20.7){\small $\B_n$}
\put(3,21.1){\small $4$}
\put(13,20.8){\small $n \ge 2$}
\put(0.5,16.7){\small $\D_n$}
\put(13,16.8){\small $n\ge 4$}
\put(0.5,12.8){\small $\E_6$}
\put(12.5,12.8){\small $\E_7$}
\put(0.5,7.9){\small $\E_8$}
\put(16.5,7.9){\small $\F_4$}
\put(0.5,3){\small $\H_3$}
\put(8.5,3){\small $\H_4$}
\put(1,0.1){\small $\I_2(p)$}
\put(4,0.5){\small $p$}
\put(20.7,8.3){\small $4$}
\put(3,3.4){\small $5$}
\put(10.8,3.4){\small $5$}
\put(6.5,0.1){\small $p \ge 5$}
\end{picture}} 

\bigskip
\centerline{{\bf Figure 2.1.} Connected spherical type Coxeter graphs.}
\end{figure}

\bigskip\noindent
Let $\Gamma$ be a Coxeter graph, and let $(W,S)$ be its Coxeter system.
Denote by $S^*$ the free monoid on $S$.
Let $w \in W$.
A word $\mu=s_1 \cdots s_\ell \in S^*$ is an {\it expression} of $w$ if the equality $w=s_1 \cdots s_\ell$ holds in $W$.
The {\it length} of $w$, denoted by $\lg(w)$, is defined to be the minimal length of an expression of $w$.
An expression $\mu=s_1 \cdots s_\ell$ of $w$ is said to be {\it reduced} if $\ell=\lg(w)$.

\bigskip\noindent
Let  $\mu,\mu'\in S^*$.
We say that there is an {\it elementary M-transformation} joining $\mu$ to $\mu'$ if there exist $\nu_1,\nu_2 \in S^*$ and $s,t \in S$ such that $m_{s,t} \neq \infty$,
\[
\mu=\nu_1\, \Pi(s,t:m_{s,t})\, \nu_2\,,\quad \text{and}\quad \mu'= \nu_1\, \Pi(t,s :m_{s,t})\, \nu_2\,.
\]

\bigskip\noindent
{\bf Theorem 2.7}
(Tits \cite{Tits2}).
{\it Let $w \in W$, and let $\mu,\mu'$ be two reduced expressions of $w$.
Then there is a finite sequence of elementary M-transformations joining $\mu$ to $\mu'$.}

\bigskip\noindent
Let $(A,\Sigma)$ be the Artin system of $\Gamma$.
Recall the epimorphism $\theta : A \to W$ which sends $\sigma_s$ to $s$ for all $s \in S$.
We define a set-section $\tau: W \to A$ of $\theta$ as follows.
Let $w \in W$.
We choose a reduced expression $\mu=s_1 \cdots s_\ell$ of $w$ and we set 
\[
\tau(w) = \sigma_{s_1} \cdots \sigma_{s_\ell}\,.
\]
By Theorem 2.7 the definition of $\tau(w)$ does not depend on the choice of the reduced expression. 
Attention: $\tau$ is a set-section.
It is not a homomorphism. 
However, it is an important tool in the study of Artin groups.

\bigskip\noindent
The following theorem is fundamental in the combinatorial study of Coxeter groups. 

\bigskip\noindent
{\bf Theorem 2.8}
(Bourbaki \cite{Bourb1}).
{\it Let $(W,S)$ be a Coxeter system.
\begin{itemize}
\item[(1)]
Let $w \in W$, let $s \in S$, and let $\mu=s_1 \cdots s_\ell$ be a reduced expression of $w$.
Then, either $\lg(ws) = \lg(w) +1$, or there exists an index $i \in \{1, \dots, \ell\}$ such that $w=s_1 \cdots \widehat{s_i} \cdots s_\ell s$.
\item[(2)]
Let $w \in W$ and $s,t \in S$.
If $\lg(ws) = \lg(tw) = \lg(w)+1$ and $\lg(tws) < \lg(ws)$, then $ws=tw$.
\end{itemize}}

\bigskip\noindent
Recall that, for $X\subset S$, we set $M_X = (m_{s,t})_{s,t \in X}$, where $M=(m_{s,t})_{s,t \in S}$ is the Coxeter matrix of $\Gamma$, we denote by $\Gamma_X$ the Coxeter graph of $M_X$, and we denote by $W_X$ the subgroup of $W=W_\Gamma$ generated by $X$.
Recall also that, by \cite{Bourb1}, the pair $(W_X,X)$ is the Coxeter system of $\Gamma_X$, and $W_X$ is called {\it standard parabolic subgroup} of $W$.
Let $X,Y$ be two subsets of $S$.
We say that an element $w \in W$ is {\it $(X,Y)$-minimal} if it is of minimal length in the double-coset $W_XwW_Y$.

\bigskip\noindent
{\bf Proposition 2.9}
(Bourbaki \cite{Bourb1}). 
{\it Let $(W,S)$ be a Coxeter system.
\begin{itemize}
\item[(1)]
Let $X,Y$ be two subsets of $S$, and let $w \in W$. 
Then there exists a unique $(X,Y)$-minimal element lying in $W_X w W_Y$.
\item[(2)]
Let $X \subset S$, and let $w \in W$. 
Then $w$ is $(\emptyset,X)$-minimal if and only if $\lg (ws) > \lg (w)$ for all $s \in X$, and $\lg (ws) > \lg (w)$ for all $s \in X$ if and only if $\lg (wu) = \lg(w) + \lg(u)$ for all $u \in W_X$.
\item[(3)]
Let $X \subset S$, and let $w \in W$. 
Then $w$ is $(X,\emptyset)$-minimal if and only if $\lg (sw) > \lg (w)$ for all $s \in X$, and $\lg (sw) > \lg (w)$ for all $s \in X$ if and only if $\lg (uw) = \lg(u) + \lg(w)$ for all $u \in W_X$.
\item[(4)]
Let $X \subset S$, and let $w \in W_X$.
If $\mu=s_1 \cdots s_\ell$ is a reduced expression of $w$, then $s_1, \dots, s_\ell \in X$.
\end{itemize}}

\subsection{Preliminaries on Vinberg systems}

In this subsection we present the main tool that we will use to pass from the Salvetti complex of a Coxeter arrangement, denoted by $\Sal(\AA)$, to the Salvetti complex of the associated Coxeter graph, denoted by $\Sal(\Gamma)$ (see Subsection 3.2).
This tool basically says that two posets are isomorphic (see Theorem 2.10 below).
The first poset is the poset of facets of the Coxeter arrangement, while the second poset, denoted by $\PP^f$, is made of the cosets of finite parabolic subgroups in $W$.
Theorem 2.10 is essentially due to Vinberg \cite{Vinbe1}, but the proofs of \cite[Chap. V]{Bourb1} can be easily adapted to prove the theorem. 

\bigskip\noindent
Let $\AA$ be a hyperplane arrangement in an nonempty open convex cone $I$ in $V= \R^\ell$.
A {\it chamber} of $\AA$ is defined to be a connected component of $I \setminus (\cup_{H \in \AA} H)$.
We denote by $\CC (\AA)$ the set of chambers of $\AA$.
For $H \in \AA$, we set $I^H = I \cap H$ and $\AA^H = \{ H' \cap H \mid  H' \in \AA \setminus \{H\} \text{ and } H' \cap H \cap I \neq \emptyset \}$.
Observe that $I^H$ is a nonempty open convex cone in $H$, and $\AA^H$ is a hyperplane arrangement in $I^H$.
For $H \in \AA$, a chamber of $\AA^H$ is called a {\it face} of $\AA$ (or {\it $1$-codimensional facet} of $\AA$).
For $d \in \N$, we define a {\it $d$-codimensional facet} of $\AA$ by induction on $d$ as follows.
The chambers of $\AA$ are the $0$-codimensional facets.
The faces of $\AA$ are the $1$-codimensional facets.
For $d \ge 2$, a $d$-codimensional facet of $\AA$ is a $(d-1)$-codimensional facet of some $\AA^H$, where $H \in\AA$.
We denote by $\FF(\AA)$ the set of all facets of $\AA$.
Observe that $\FF(\AA)$ is a partition of $I$.
For $F \in \FF(\AA)$, we denote by $\bar F$ the closure of $F$ in $I$.
Then $\FF(\AA)$ is endowed with the partial order relation $\preceq$ defined by $F_1 \preceq F_2$ if $F_1 \subseteq \bar F_2$.

\bigskip\noindent
Let $F$ be a $d$-codimensional facet.
Define the {\it support} of $F$, denoted by $|F|$, to be the linear subspace of $V$ spanned by $F$.
Set $I^F = |F| \cap I$, $\AA_F = \{ H \in \AA \mid H \supset F\}$, and $\AA^F = \{ H \cap |F| \mid H \in \AA \setminus \AA_F \text{ and }H \cap I^F \neq \emptyset\}$.
Observe that $|F|$ is a $d$-codimensional linear subspace of $V$, $I^F$ is a nonempty open convex cone in $|F|$, $\AA^F$ is a hyperplane arrangement in $I^F$, and $F$ is a chamber of $\AA^F$.
On the other hand, $\AA_F$ is is a finite hyperplane arrangement in $I$.
Moreover, if $d \ge 1$, we have $\cap_{H \in \AA_F} H = |F|$.
For $d=0$, we set $|F|=\cap_{H \in \AA_F} H = \cap_{H \in \emptyset} H = V$.

\bigskip\noindent
{\bf Example.}
Set $V=\R^3$ and $I=\{(x,y,z) \in V \mid z>0\}$.
For $k \in \Z$, denote by $H_k$ the plane of $V$ of equation $x=kz$, and denote by $H_k'$ the plane of equation $y=kz$.
Set $\AA=\{ H_k,\ H_k' \mid k \in \Z\}$.
This is a hyperplane arrangement in $I$.
The trace of $\AA$ in the affine plane of equation $z=1$ is pictured in Figure 2.2.
Let $F=\{(0,0,z) \mid z>0\}$.
Then $F$ is a $2$-codimensional facet of $\AA$ whose support is the line $|F|$ of equations $x=y=0$.
Here we have $I^F=F$, $\AA^F=\emptyset$, and $\AA_F = \{H_0, H_0'\}$.
Observe that the set of facets $F' \in \FF(\AA)$ satisfying $F \preceq F'$ is made of $4$ chambers, $4$ faces, and $F$ itself.

\begin{figure}[tbh]
\bigskip
\centerline{
\setlength{\unitlength}{0.5cm}
\begin{picture}(9,9)
\put(1,1){\includegraphics[width=4cm]{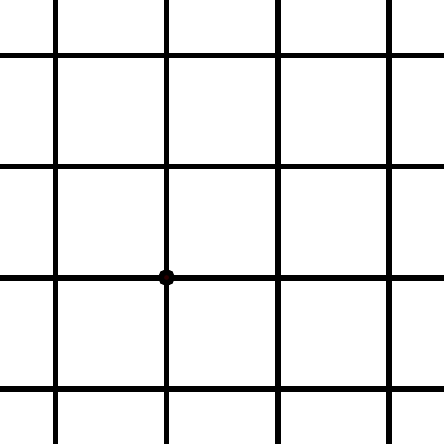}}
\put(1.4,0.2){\small $H_{-1}$}
\put(3.7,0.2){\small $H_0$}
\put(5.7,0.2){\small $H_1$}
\put(7.7,0.2){\small $H_2$}
\put(-0.5,1.8){\small $H_{-1}'$}
\put(0,3.8){\small $H_0'$}
\put(0,5.8){\small $H_1'$}
\put(0,7.8){\small $H_2'$}
\put(4.2,4.2){\small $F$}
\end{picture}} 

\bigskip
\centerline{{\bf Figure 2.2.} A facet.}
\end{figure}

\bigskip\noindent
For a given Coxeter graph $\Gamma$ and its Coxeter system $(W,S)$, we set $\SS^f = \SS^f_\Gamma =\{X \subset S \mid W_X \text{ is finite}\}$ and $\PP^f=\PP^f_\Gamma =\{ wW_X \mid w \in W \text{ and } X \in \SS^f\}$.
The set $\PP^f$ is assumed to be ordered by the inclusion.

\bigskip\noindent
{\bf Theorem 2.10}
(Vinberg \cite{Vinbe1}).
{\it Let $(W,S)$ be a Vinberg system, let $C_0$ be its fundamental chamber, let $I$ be its Tits cone, and let $\AA$ be its Coxeter arrangement. 
Denote by $\Gamma$ the Coxeter graph of $(W,S)$, viewed as a Coxeter system.}

\bigskip\noindent
{\bf I.}
{\it Let $\FF(C_0) = \{ F \in \FF(\AA) \mid F \preceq C_0 \}$.
Then there is a bijection $\iota : \SS_\Gamma^f \to \FF(C_0)$ such that 
\[
\bigcap_{s \in X} H_s = |\iota(X)|
\]
for all $X \in \SS_\Gamma^f$.
Moreover, the following properties hold.
\begin{itemize}
\item[(1)]
Let $X,Y \in \SS_\Gamma^f$.
We have $X \subset Y$ if and only if $\iota(Y) \preceq \iota(X)$.
\item[(2)]
For $X \in \SS_\Gamma^f$, the stabilizer $\{w \in W \mid w(\iota(X))=\iota(X)\}$ of $\iota(X)$ is equal to $W_X$, and every element of $W_X$ pointwise fixes $\iota(X)$.
\end{itemize}}

\bigskip\noindent
{\bf II.}
{\it There is a bijection $\tilde \iota : \PP_\Gamma^f \to \FF(\AA)$ defined by 
\[
\tilde \iota(wW_X) = w(\iota(X))\,.
\]
Moreover, the following properties hold.
\begin{itemize}
\item[(1)]
Let $u,v \in W$ and $X,Y \in \SS_\Gamma^f$.
We have $uW_X \subset vW_Y$ if and only if $\tilde \iota (vW_Y) \preceq \tilde \iota (uW_X)$.
\item[(2)]
Let $u,v \in W$ and $X,Y \in \SS_\Gamma^f$.
We have $uW_X \subset vW_Y$ if and only if $X \subset Y$ and $u \in vW_Y$.
\item[(3)]
The restriction of $\tilde \iota$ to $W$ is the bijection
\[
\begin{array}{ccc}
W & \to & \CC(\AA)\\
w & \mapsto & w(C_0)
\end{array}
\]
\end{itemize}}

\bigskip\noindent
{\bf III.}
{\it Let $X$ be a subset of $S$, and let $w$ be an element of $W$. 
Then $w$ is $(\emptyset,X)$-minimal if and only if $H_r$ does not separate $C_0$ and $w^{-1}(C_0)$ for every reflection $r$ lying in $W_X$.}

\subsection{Preliminaries on Artin monoids}

In this subsection we present some results on Artin monoids, that we will need to prove Theorem~4.10 ($K(\pi,1)$ conjecture for spherical type Artin groups).
Most of the results of the subsection come from \cite{BriSai1} (see also \cite{Miche1}), and their proofs are independent from the techniques presented here.

\bigskip\noindent
The {\it Artin monoid} of a Coxeter graph $\Gamma$ is the monoid $A_\Gamma^+$ defined by the following monoid presentation.
\[
A_\Gamma^+ = \langle \Sigma \mid \Pi(\sigma_s,\sigma_t : m_{s,t}) = \Pi(\sigma_t,\sigma_s : m_{s,t}) \text{ for all } s,t \in S,\ s \neq t,\ m_{s,t} \neq \infty \rangle^+ \,.
\]
By \cite{Paris2}, the natural homomorphism $A_\Gamma^+ \to A_\Gamma$ is injective.

\bigskip\noindent
Since the relations that define $A_\Gamma^+$ are homogeneous, $A_\Gamma^+$ is endowed with a length function $\lg: A_\Gamma^+ \to \N$ that associates to each element $\alpha \in A_\Gamma^+$ the length of any expression of $\alpha$ on the elements of $\Sigma$.
Note that $\lg(\alpha \beta) = \lg(\alpha) + \lg (\beta)$ for all $\alpha, \beta \in A_\Gamma^+$.
For $\alpha, \beta \in A_\Gamma^+$, we set $\alpha \preceq_L \beta$ if there exists $\gamma \in A_\Gamma^+$ such that $\alpha\gamma=\beta$.
Similarly, we set $\alpha\preceq_R \beta$ if there exists $\gamma \in A_\Gamma^+$ such that $\gamma \alpha = \beta$.
Note that the existence of the length function implies that $\preceq_L$ and $\preceq_R$ are partial order relations on $A_\Gamma^+$.

\bigskip\noindent
{\bf Theorem 2.11}
(Brieskorn, Saito \cite{BriSai1}).
{\it Let $\Gamma$ be a Coxeter graph, and let $E$ be a nonempty finite subset of $A_\Gamma^+$.
\begin{itemize}
\item[(1)]
$E$ has a greater lower bound for the relation $\preceq_L$ (resp. $\preceq_R$), denoted by $\wedge_LE$ (resp. $\wedge_RE$). 
\item[(2)]
If $E$ has an upper bound for the relation $\preceq_L$ (resp. $\preceq_R$), then $E$ has a least upper bound for the relation $\preceq_L$ (resp. $\preceq_R$), denoted by $\vee_LE$ (resp. $\vee_RE$).
\end{itemize}}

\bigskip\noindent
Recall the natural homomorphism $\theta : A_\Gamma \to W_\Gamma$ and its set-section $\tau : W_\Gamma \to A_\Gamma$.
Notice that the image of $\tau$ is contained in $A_\Gamma^+$.
We complete the above theorem with the following. 

\bigskip\noindent
{\bf Proposition 2.12}
(Brieskorn, Saito \cite{BriSai1}).
{\it Let $\Gamma$ be a Coxeter graph.
\begin{itemize}
\item[(1)]
Let $E$ be a nonempty finite subset of $\tau(W_\Gamma)$. 
If $\vee_LE$ (resp. $\vee_RE$) exists, then $\vee_LE \in \tau(W_\Gamma)$ (resp. $\vee_RE \in \tau(W_\Gamma)$).
\item[(2)]
Let $X$ be a subset of $S$.
Recall that $\Sigma_X$ denotes the set $\Sigma_X=\{ \sigma_s \mid s \in X\}$.
Then $\vee_L\Sigma_X$  (resp. $\vee_R \Sigma_X$) exists if and only if $W_X$ is finite (that is, $X \in \SS^f$).
\end{itemize}}

\bigskip\noindent
The last preliminary result on Artin monoids that we will need concerns only the spherical type ones. 

\bigskip\noindent
{\bf Theorem 2.13}
(Brieskorn, Saito \cite{BriSai1}, Deligne \cite{Delig1}).
{\it Let $\Gamma$ be a spherical type Coxeter graph.
Set $\Delta=\vee_L\Sigma$ (this element exists by Proposition 2.12).
Then $\Delta=\vee_R\Sigma$, and every element $\beta\in A_\Gamma$ can be written in the form $\beta = \Delta^{-k} \alpha$ with $\alpha \in A_\Gamma^+$ and $k \in \N$.}

%%%%%%%%%%%%%%
%%CHAPITRE 3%%
%%%%%%%%%%%%%%

\section{Salvetti complexes}

\subsection{Salvetti complex of a hyperplane arrangement}

In this subsection $I$ denotes an nonempty open convex cone in a real vector space $V$ of dimension $\ell$, and $\AA$ denotes a hyperplane arrangement in $I$.
Our aim is to define a (geometric) simplicial complex $\Sal(\AA)$ and to prove that $\Sal(\AA)$ has the same homotopy type as $M(\AA)$. 
We start recalling some definitions from the previous section. 

\bigskip\noindent
The arrangement $\AA$ determines a partition of $I$ into facets. 
We denote by $\FF(\AA)$ the set of facets, and by $\CC(\AA)$ the set of chambers ($0$-codimensional facets) of $\AA$. 
We order $\FF(\AA)$ by $F_1 \preceq F_2$ if $F_1 \subseteq \bar F_2$,
where, for $F \in \FF(\AA)$, $\bar F$ denotes the closure of $F$ in $I$. 
The support of a facet $F$, denoted by $|F|$, is the linear subspace of $V$ spanned by $F$.
We set $I^F= I \cap |F|$, $\AA_F=\{H \in \AA \mid F \subset H\}$, and $\AA^F=\{ H \cap |F| \mid H \in \AA \setminus \AA_F \text{ and } H \cap I^F \neq \emptyset\}$.
Finally, for $F \in \FF(\AA)$ and $C \in \CC(\AA)$, we denote by $C_F$ the chamber of $\AA_F$ containing $C$.

\bigskip\noindent
{\bf Example.}
We go back to the example of the previous section.
We set $V=\R^3$ and $I=\{(x,y,z) \in V \mid z>0\}$.
For $k \in \Z$, we denote by $H_k$ the plane of $V$ of equation $x=kz$, and by $H_k'$ the plane of equation $y=kz$, and we set $\AA=\{ H_k,\ H_k' \mid k \in \Z\}$.
Consider the facet $F=\{(0,0,z) \mid z>0\}$.
Let $C$ be the cone over a square bounded by $H_0,H_1, H_0',H_1'$.
Then $C$ is a chamber of $\AA$, we have $F \preceq C$, and $C_F$ is the cone $\{(x,y,z) \in V \mid x>0,\ y>0,\ z>0\}$ (see Figure 3.1).

\begin{figure}[tbh]
\bigskip
\centerline{
\setlength{\unitlength}{0.5cm}
\begin{picture}(21,9)
\put(1,1){\includegraphics[width=10cm]{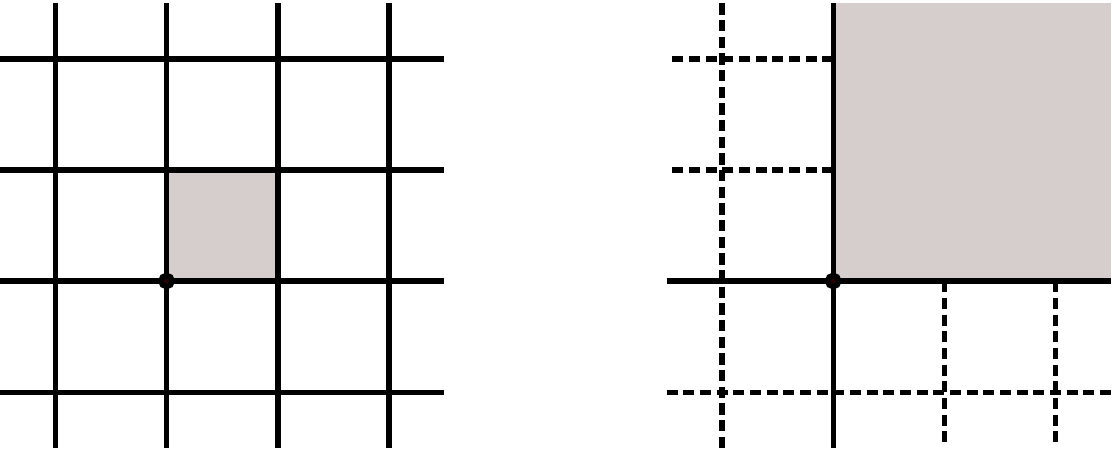}}
\put(1.4,0.2){\small $H_{-1}$}
\put(3.7,0.2){\small $H_0$}
\put(5.7,0.2){\small $H_1$}
\put(7.7,0.2){\small $H_2$}
\put(-0.5,1.8){\small $H_{-1}'$}
\put(0,3.8){\small $H_0'$}
\put(0,5.8){\small $H_1'$}
\put(0,7.8){\small $H_2'$}
\put(3.2,3.2){\small $F$}
\put(4.8,4.8){\small $C$}
\put(15.7,0.2){\small $H_0$}
\put(12,3.8){\small $H_0'$}
\put(18.2,6){\small $C_F$}
\end{picture}} 

\bigskip
\centerline{{\bf Figure 3.1.} $C_F$: an example.}
\end{figure}

\bigskip\noindent
We set
\[
\Sal_0(\AA) = \{(F,C) \in \FF(\AA) \times \CC(\AA) \mid F \preceq C\}\,.
\]
We define a relation $\preceq$ on $\Sal_0(\AA)$ as follows.
\begin{center}
$(F,C) \preceq (F',C')$ if $F \preceq F'$ and $C_F \subset C'_{F'}$.
\end{center}
It is easily checked that $\preceq$ is an order relation on $\Sal_0(\AA)$.
We define the {\it Salvetti complex} of $\AA$, denoted by $\Sal(\AA)$, as the geometric realization of the derived complex of $(\Sal_0(\AA), \preceq)$.

\bigskip\noindent
{\bf Remark.}
If $\AA$ is the Coxeter arrangement of a Vinberg system $(W,S)$, then $W$ acts on $\Sal_0(\AA)$ as follows.
\[
w\,(F,C) = (wF,wC)\,,
\]
for $w \in W$ and $(F,C) \in \Sal_0(\AA)$.
The ordering $\preceq$ is invariant under the action of $W$, thus this action induces an action of $W$ on $\Sal(\AA)$.

\bigskip\noindent
{\bf Theorem 3.1.}
{\it There exists a homotopy equivalence $f: \Sal(\AA) \to M(\AA)$.
Moreover, if $\AA$ is the Coxeter arrangement of a Vinberg system $(W,S)$, then $f$ is equivariant under the actions of $W$ and induces a homotopy equivalence $\bar f : \Sal(\AA)/W \to M(\AA)/W$.}

\bigskip\noindent
{\bf Proof.}
We shall define a family $\{U(F,C) \mid (F,C) \in \Sal_0(\AA) \}$ of open subsets of $M(\AA)$, and we shall prove the following. 
\begin{itemize}
\item[(1)]
Let $(F,C), (G,D) \in \Sal_0(\AA)$. 
If $U(F,C) = U(G,D)$, then $(F,C) = (G,D)$.
\item[(2)]
We have 
\[
M(\AA) = \bigcup_{(F,C) \in \Sal_0(\AA)} U(F,C)\,.
\]
\item[(3)]
Let $(F_0,C_0), (F_1,C_1), \dots, (F_p,C_p) \in \Sal_0(\AA)$.
We have 
\[
U(F_0,C_0) \cap U(F_1, C_1) \cap \cdots \cap U(F_p,C_p) \neq \emptyset
\]
if and only if, up to permutation, we have a chain
\[
(F_0,C_0) \prec (F_1,C_1) \prec \cdots \prec (F_p,C_p)\,.
\]
\item[(4)]
Let $(F_0,C_0) \prec \cdots \prec (F_p,C_p)$ be a chain in $\Sal_0(\AA)$. 
Then $U(F_0,C_0) \cap U(F_1,C_1) \cap \cdots \cap U(F_p,C_p)$ is contractible.
\end{itemize}
Moreover, if $\AA$ is the Coxeter arrangement of a Vinberg system $(W,S)$, we shall prove the following. 
\begin{itemize}
\item[(5)]
Let $(F,C) \in \Sal_0(\AA)$ and $w \in W\setminus \{1\}$. 
Then $w\,U(F,C)=U(w\,F, w\,C)$ and $U(F,C) \cap U(w\,F,w\,C)=\emptyset$.
\end{itemize}
By Theorem 2.1 and Proposition 2.2, Theorem 3.1 will be a straightforward consequence of (1)--(5). 

\bigskip\noindent
For $(F,C) \in \Sal_0(\AA)$, the open subset $U(F,C)$ will be of the form $U(F,C) = \omega(F) \times C_F$, where $\omega (F)$ is an open subset of $I$.
We turn now to construct $\omega(F)$ and study its properties.

\bigskip\noindent
Throughout the proof we adopt the following definitions and notations. 
A {\it chain of length} $p+1$ in $\FF(\AA)$ is a sequence $(F_0,F_1, \dots, F_p)$ in $\FF(\AA)$ such that $F_0 \prec F_1 \prec \cdots \prec F_p$.
We set $\gamma \le \gamma'$ if $\gamma= (F_0, F_1, \dots, F_p)$ and $\gamma' = (F_0',F_1', \dots, F_q')$ are two chains in $\FF(\AA)$ such that $F_0 = F_0'$ and $\{F_1, \dots, F_p\} \subseteq \{F_1', \dots, F_q'\}$.
For $F \in \FF(\AA)$, we denote by $\Chain(F)$ the set of chains $\gamma=(F_0,F_1, \dots, F_p)$ such that $F_0=F$.
More generally, if $\gamma$ is a chain, we denote by $\Chain(\gamma)$ the set of chains $\gamma'$ such that $\gamma \le \gamma'$.

\bigskip\noindent
For all $F \in \FF(\AA)$, we fix a point $x(F) \in F$.
If $\AA$ is the Coxeter arrangement of a Vinberg system $(W,S)$, we choose the points $x(F)$ so that  $w\,x(F)=x(w\,F)$ for all $F \in \FF(\AA)$ and all $w \in W$.
For a given chain $\gamma =(F_0,F_1, \dots, F_p)$ in $\FF(\AA)$ we set 
\[
\Delta(\gamma) = \{y + t_1x(F_1) + \cdots + t_px(F_p) \mid y\in F_0 \text{ and }t_1, \dots, t_p >0\}\,.
\]
Note that $\Delta(F_0,F_1, \dots, F_p) \subset F_p$.
Note also that, if $\AA$ is the Coxeter arrangement of a Vinberg system $(W,S)$, then
\[
w\, \Delta(F_0,F_1, \dots, F_p) = \Delta(w\,F_0, w\,F_1, \dots, w\, F_p)
\]
for every chain $(F_0,F_1, \dots, F_p)$ and every $w \in W$.

\bigskip\noindent
For a given $F \in \FF(\AA)$ we set 
\[
\omega(F) = \bigcup_{\gamma \in \Chain(F)} \Delta(\gamma)
\]
(see Figure 3.2). More generally, for a given chain $\gamma$ in $\FF(\AA)$, we set
\[
\omega(\gamma) = \bigcup_{\gamma' \in \Chain(\gamma)} \Delta(\gamma')\,.
\]

\begin{figure}[tbh]
\bigskip
\centerline{
\setlength{\unitlength}{0.5cm}
\begin{picture}(10,10)
\put(0,0){\includegraphics[width=5cm]{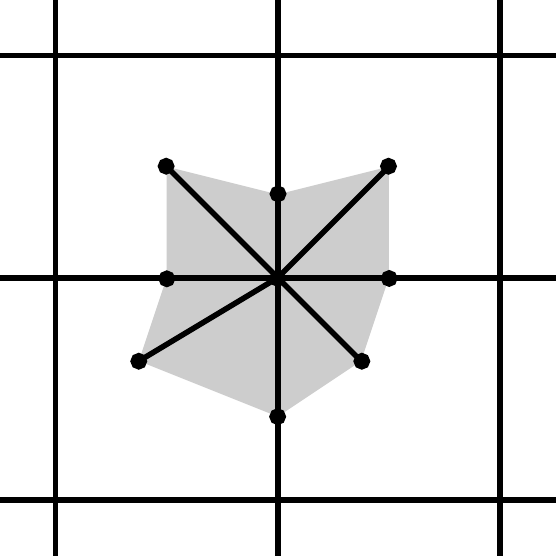}}
\put(4.3,4){\small $F$}
\end{picture}} 

\bigskip
\centerline{{\bf Figure 3.2.} The set $\omega(F)$.}
\end{figure}

\bigskip\noindent
{\bf Claim 1.}
{\it Let $F \in \FF(\AA)$.
Then $\omega(F)$ is an open subset of $I$.}

\bigskip\noindent
{\it Proof of Claim 1.}
Let $d$ be the codimension of $F$.
For $k \ge 0$, we denote by $I_{d-k}$ the union of facets of $\FF(\AA)$ of codimension $\ge d-k$, and we prove by induction on $k$ that $I_{d-k} \cap \omega(F)$ is an open subset of $I_{d-k}$.
The set $I_d \cap \omega(F)=F$ is obviously open in $I_d$, thus we may assume that $k \ge 1$ plus the induction hypothesis.
If $X$ is a subset of the cone $I$ and $x$ is a point in $I$, the following set
\[
\{y+tx \mid y \in X \text{ and } t>0 \}
\]
is called the {\it open cone over $X$ with direction $x$}.
Let $G$ be a facet of codimension $d-k$.
We denote by $\partial G$ the union of the facets $K$ such that $K \prec G$.
If $F \not\prec G$, then $G \cap \omega(F) = \emptyset$.
If $F \prec G$, then $\partial G \cap \omega(F)$ is an open subset of $\partial G$ (by induction) and $G \cap \omega(F)$ is the open cone over $\partial G \cap \omega(F)$ with direction $x(G)$.
This implies that $I_{d-k} \cap \omega(F)$ is an open subset of $I_{d-k}$.

\bigskip\noindent
{\bf Claim 2.}
{\it Let $\gamma=(F_0,F_1, \dots, F_p)$ and $\gamma'= (F_0', F_1', \dots, F_q')$ be two chains in $\FF(\AA)$.
If $q\le p$ and $\Delta(\gamma) \cap \Delta(\gamma') \neq \emptyset$, then $(F_0', \dots, F_q') = (F_{p-q}, \dots, F_p)$.}

\bigskip\noindent
{\it Proof of Claim 2.}
We argue by induction on $q$.
Suppose $q=0$.
Let $z \in \Delta(F_0, F_1, \dots, F_p) \cap \Delta(F_0')$.
Then $z \in F_p \cap F_0'$, thus $F_p \cap F_0' \neq \emptyset$, therefore $F_p=F_0'$.

\bigskip\noindent
Suppose $q >0$ plus the induction hypothesis.
Let $z \in \Delta(F_0, \dots, F_p) \cap \Delta(F_0', \dots, F_q')$.
Notice that $z \in F_p \cap F_q'$, thus $F_p \cap F_q' \neq \emptyset$, hence $F_p=F_q'$.
We write $z$ in the form $z=y + t_1\,x(F_1) + \cdots + t_p\,x(F_p)$, where $y \in F_0$ and $t_1, \dots, t_p >0$.
Similarly, we write $z = y' + t_1'\, x(F_1') + \cdots + t_q'\, x(F_q')$, where $y' \in F_0'$ and $t_1', \dots, t_q' >0$.
Let $z_1 =y + t_1\,x(F_1) + \cdots + t_{p-1}\,x(F_{p-1})$, and let $z_1' = y' + t_1'\, x(F_1') + \cdots + t_{q-1}'\, x(F_{q-1}')$.
Notice that $z_1 \in F_{p-1}$ and $z_1' \in F_{q-1}'$.
If $t_p >t_q'$, then 
\[
z_1' = y + t_1\,x(F_1) + \cdots + t_{p-1}\,x(F_{p-1}) + (t_p-t_q')x(F_p) \in F_p=F_q'\,.
\]
This is a contradiction since $z_1' \in F_{q-1}'$.
We prove in the same way that the inequality $t_p<t_q'$ cannot hold.
It follows that $t_p=t_q'$, thus $z_1=z_1'$.
By induction, we conclude that $(F_0', \dots, F_{q-1}') = (F_{p-q}, \dots, F_{p-1})$.

\bigskip\noindent
{\bf Claim 3.}
{\it Let $F \in \FF(\AA)$, and let $z \in \omega(F)$. 
There exists a unique chain $\gamma \in \Chain(F)$ such that $z \in \Delta(\gamma)$.}

\bigskip\noindent
{\it Proof of Claim 3.}
Let $\gamma = (F_0, \dots, F_p)$ and $\gamma' = (F_0', \dots, F_q')$ be elements of $\Chain(F)$ such that $z \in \Delta(\gamma) \cap \Delta(\gamma')$.
We can assume without loss of generality that $q \le p$.
By Claim 2, we have $(F_0', \dots, F_q') = (F_{p-q}, \dots, F_p)$.
Since, moreover, $F_0=F_0'=F$, it follows that $p=q$ and $\gamma=\gamma'$.

\bigskip\noindent
{\bf Claim 4.}
{\it Let $F,G \in \FF(\AA)$.
If $\omega(F) \cap \omega(G) \neq \emptyset$, then either $F \preceq G$, or $G \preceq F$.}

\bigskip\noindent
{\it Proof of Claim 4.}
Let $F,G \in \FF(\AA)$ such that $\omega(F) \cap \omega(G) \neq \emptyset$.
There exist a chain $\gamma = (F_0,F_1, \dots, F_p)$ lying in $\Chain(F)$ and a chain $\gamma' = (F_0', \dots, F_q')$ lying in $\Chain(G)$ such that $\Delta(\gamma) \cap \Delta(\gamma') \neq \emptyset$.
We can assume without loss of generality that $q \le p$.
By Claim 2, we have $(F_0', \dots, F_q') = (F_{p-q}, \dots, F_p)$.
Hence, $F=F_0 \preceq F_{p-q}=F_0'=G$. 

\bigskip\noindent
A straightforward consequence of Claim 4 is the following.

\bigskip\noindent
{\bf Claim 5.}
{\it Let $F_0,F_1, \dots, F_p \in \FF(\AA)$.
If $\omega(F_0) \cap \cdots \cap \omega(F_p) \neq \emptyset$, then, up to a permutation of the indices, we have a chain $F_0 \prec F_1 \prec \cdots \prec F_p$.}

\bigskip\noindent
{\bf Claim 6.}
{\it Let $\gamma = (F_0, \dots, F_p)$ be a chain in $\FF(\AA)$.
Then
\[
\omega(F_0) \cap \cdots \cap \omega(F_p) = \omega(\gamma)\,.
\]}

\bigskip\noindent
{\it Proof of Claim 6.}
Since the inclusion $\omega(\gamma) \subset \omega(F_0) \cap \cdots \cap \omega(F_p)$ is obvious, we only need to prove $\omega(F_0) \cap \cdots \cap \omega(F_p) \subset \omega(\gamma)$. 
Let $z \in \omega(F_0) \cap \cdots \cap \omega(F_p)$.
Since $z \in \omega(F_0)$, there exists a chain $\delta = (G_0, \dots, G_q)$ lying in $\Chain(F_0)$ such that $z \in \Delta(G_0, \dots, G_q)$.
Let $i \in \{1, \dots, p\}$.
There is also a chain $\delta'=(G_0', \dots, G_r')$ lying in $\Chain(F_i)$ such that $z \in \Delta(\delta')$.
If $q<r$, then, by Claim 2, $F_0=G_0 = G_{r-q}'$.
But, this is not possible because, otherwise, we would have $F_i = G_0' \preceq G_{r-q}'=F_0$ and $F_0 \prec F_i$.
So, $r\le q$ and, again by Claim 2, $F_i = G_0' = G_{q-r}$.
This shows that $\gamma \le \delta$, thus $z \in \Delta(\delta) \subset \omega(\gamma)$. 

\bigskip\noindent
{\bf Claim 7.}
{\it Let $\gamma = (F_0, \dots, F_p)$ be a chain in $\FF(\AA)$.
Then $\omega(\gamma)$ is contractible.}

\bigskip\noindent
{\it Proof of Claim 7.}
We choose $y_0 \in F_0$ and we set 
\[
z_0 = y_0 + x(F_1) + \cdots + x(F_p)\,.
\]
For $t \in [0,1]$, we define the map $h_t : \omega(\gamma) \to \omega(\gamma)$ as follows.
Let $z \in \omega(\gamma)$.
Let $\delta = (G_0, \dots, G_q) \in \Chain(\gamma)$ such that $z \in \Delta(\delta)$.
There exist $y \in F_0=G_0$ and $t_1, \dots, t_q >0$ such that $z=y + t_1\,x(G_1) + \cdots + t_q\, x(G_q)$.
We set
\[
h_t(z) = t\,y+(1-t)y_0 + u_1(t)\, x(G_1) + \cdots + u_q(t)\,x(G_q)\,,
\]
where
\[
u_j(t) = \left\{ \begin{array}{ll}
t\,t_j & \text{if } G_j \not\in \{F_1, \dots, F_p\}\\
t\,t_j + (1-t) & \text{if } G_j \in \{F_1, \dots, F_p\}
\end{array} \right.
\]
It is easily seen that the map
\[
\begin{array}{ccc}
\omega(\gamma) \times [0,1] & \to & \omega(\gamma)\\
(z,t) & \mapsto & h_t(z)
\end{array}
\]
is well-defined and continuous.
Moreover, we have $h_t(z_0)=z_0$ for all $t \in [0,1]$, $h_1=\Id$, and $h_0(z) = z_0$ for all $z \in \omega(\gamma)$.
This shows that $\{z_0\}$ is a deformation retract of $\omega(\gamma)$.

\bigskip\noindent
{\bf Claim 8.}
{\it Suppose that $\AA$ is the Coxeter arrangement of a Vinberg system $(W,S)$.
For $F \in \FF(\AA)$, we denote by $W_F=\{w \in W \mid w\,F=F\}$ the stabilizer of $F$.
Let $F \in \FF(\AA)$, and let $w \in W$.
If $\omega(F) \cap w\, \omega(F) \neq \emptyset$, then $w \in W_F$.}

\bigskip\noindent
{\it Proof of Claim 8.}
Suppose that $\omega(F) \cap w\, \omega(F) \neq \emptyset$.
Let $z \in \omega(F) \cap w\, \omega(F)$.
There are chains $\gamma = (F_0, \dots, F_p)$ and $\gamma' = (F_0', \dots, F_q')$ lying in $\Chain(F)$ such that
\[
z \in \Delta(\gamma) \cap w\, \Delta(\gamma') = \Delta(F_0, \dots, F_p) \cap \Delta (w\, F_0', \dots w\, F_q')\,.
\]
Assume $q \le p$.
The case $p \le q$ is proved in the same way.
By Claim 2, $(w\, F_0', \dots, w\, F_q') = (F_{p-q}, \dots, F_p)$.
Since $\codim\, F_0 =\codim\,w\,F_0' = \codim\,F$, it follows that $p=q$ and $w\,\gamma' = \gamma$.
In particular, $F = F_0 = w\,F_0' = w\,F$, thus $w \in W_F$.

\bigskip\noindent
For $(F,C) \in \Sal_0(\AA)$, we set
\[
U(F,C) = \omega(F) \times C_F\,.
\]
We turn now to prove in the following claims that the set $\{ U(F,C) \mid (F,C) \in \Sal_0(\AA) \}$ satisfies (1)--(5). 

\bigskip\noindent
{\bf Claim 9.}
{\it Let $(F,C) \in \Sal_0(\AA)$.
Then $U(F,C) \subset M(\AA)$.}

\bigskip\noindent
{\it Proof of Claim 9.}
Let $(x,y) \in U(F,C)$.
Let $G \in \FF(\AA)$ such that $x \in G$.
Since $x \in \omega(F)$, we have $F \preceq G$.
Let $H \in \AA$.
If $H \in \AA_F$, then $y \not\in H$, since $y \in C_F$, thus $(x,y) \not\in H \times H$.
If $H \not\in \AA_F$, then $H \not\in \AA_G$, since $\AA_G \subset \AA_F$.
But, $\AA_G=\{H' \in \AA \mid x \in H'\}$, thus $x \not\in H$, therefore $(x,y) \not\in H \times H$.
This shows that $(x,y) \in M(\AA)$.

\bigskip\noindent
{\bf Claim 10.}
{\it Let $(F,C), (G,D) \in \Sal_0(\AA)$. 
If $U(F,C) = U(G,D)$, then $(F,C)=(G,D)$.}

\bigskip\noindent
{\it Proof of Claim 10. }
We have $F \subset \omega(F)$ and, if $F' \in \FF(\AA)$ intersects $\omega(F)$, then $F \preceq F'$, thus $\omega(F)$ determines $F$.
This implies that $F=G$.
We have $\omega(F) \cap C_F \subset C$, thus $U(F,C)=\omega(F) \times C_F$ determines $C$.
This implies that $C=D$.

\bigskip\noindent
{\bf Claim 11.}
{\it We have
\[
M(\AA) \subset \bigcup_{(F,C) \in \Sal_0(\AA)} U(F,C)\,.
\]}

\bigskip\noindent
{\it Proof of Claim 11.}
Let $(x,y) \in M(\AA)$.
Let $F \in \FF(\AA)$ such that $x \in F$.
Let $H \in \AA_F$.
Since $x \in H$, we have $y \not\in H$.
Hence, there exists a chamber $\tilde C$ of $\AA_F$ such that $y \in \tilde C$.
Let $C \in \CC(\AA)$ such that $F \preceq C$ and $C_F= \tilde C$.
Then
\[
(x,y) \in F \times C_F \subset \omega(F) \times C_F = U(F,C)\,.
\]

\bigskip\noindent
{\bf Claim 12.}
{\it Let $(F,C), (G,D) \in \Sal_0(\AA)$.
If $U(F,C) \cap U(G,D) \neq \emptyset$, then either $(F,C) \preceq (G,D)$, or $(G,D) \preceq (F,C)$.}

\bigskip\noindent
{\it Proof of Claim 12.}
Suppose that $U(F,C) \cap U(G,D) \neq \emptyset$.
We have $\omega(F) \cap \omega(G) \neq \emptyset$, thus, by Claim 4, either $F \preceq G$, or $G \preceq F$.
We can assume without loss of generality that $F \preceq G$.
Then $\emptyset \neq C_F \cap D_G \subset C_G \cap D_G$, thus $C_G=D_G$, therefore $C_F \subset D_G$.

\bigskip\noindent
A straightforward consequence of Claim 12 is the following.

\bigskip\noindent
{\bf Claim 13.}
{\it Let $(F_0,C_0),(F_1,C_1), \dots, (F_p,C_p) \in \Sal_0(\AA)$.
If $U(F_0,C_0) \cap \cdots \cap U(F_p,C_p) \neq \emptyset$, then, up to a permutation of the indices, we have a chain 
\[
(F_0,C_0) \preceq (F_1,C_1) \preceq \cdots \preceq (F_p,C_p)\,.
\]}

\bigskip\noindent
{\bf Claim 14.}
{\it Let $((F_0,C_0), \dots, (F_p,C_p))$ be a chain in $\Sal_0(\AA)$.
Then
\[
U(F_0,C_0) \cap U(F_1,C_1) \cap \cdots \cap U(F_p,C_p) = \omega(F_0, \dots, F_p) \times (C_0)_{F_0}\,.
\]}

\bigskip\noindent
{\it Proof of Claim 14.}
Since $(C_0)_{F_0} \subset (C_i)_{F_i}$ for all $i \in \{1, \dots, p\}$, by Claim 6 we have
\[
U(F_0,C_0) \cap \cdots \cap U(F_p,C_p) = (\omega(F_0) \cap \cdots \cap \omega(F_p)) \times (C_0)_{F_0} = \omega(F_0, \dots, F_p) \times (C_0)_{F_0}\,.
\]

\bigskip\noindent
{\bf Claim 15.}
{\it Let $((F_0,C_0), \dots, (F_p,C_p))$ be a chain in $\Sal_0(\AA)$.
Then $U(F_0,C_0) \cap \cdots \cap U(F_p,C_p)$ is nonempty and contractible.}

\bigskip\noindent
{\it Proof of Claim 15.}
$\omega(F_0, \dots, F_p)$ is contractible by Claim 7, and $C_F$ is contractible since it is convex. 
Both spaces are obviously nonempty.

\bigskip\noindent
{\bf Claim 16.}
{\it Suppose $\AA$ is the Coxeter arrangement of a Vinberg system $(W,S)$.
Let $(F,C) \in \Sal_0(\AA)$, and let $w \in W \setminus \{1\}$.
Then $w\, U(F,C) \in \Sal_0(\AA)$ and $U(F,C) \cap w\, U(F,C) =\emptyset$.}

\bigskip\noindent
{\it Proof of Claim 16.}
We have $w\, U(F,C) = U(w\, F, w\,C) \in \Sal_0(\AA)$.
Recall that $W_F=\{ w \in W \mid wF=F\}$.
Since $W$ acts freely on $\CC(\AA)$, the group $W_F$ acts freely on $\{C \in \CC(\AA) \mid F \preceq C\}$.
If $w \not\in W_F$, then, by Claim 8, $\omega(F) \cap w\,\omega(F) = \emptyset$. 
If $w \in W_F$, then $F \preceq w\,C$ and $C \neq w\,C$, thus $C_F \neq w\,C_F$, therefore $C_F \cap w\,C_F = \emptyset$.
In both cases we have $U(F,C) \cap U(w\, F,w\, C) = \emptyset$.
\qed

\subsection{Salvetti complex of a Coxeter system}

Recall that, for a given Coxeter system $(W,S)$, we denote by $\SS^f$ the set of subsets $X$ of $S$ such that $W_X$ is finite.
The following lemma is a preliminary to the definition of the Salvetti complex of a Coxeter graph.  

\bigskip\noindent
{\bf Lemma 3.2.}
{\it Let $\Gamma$ be a Coxeter graph, and let $(W,S)$ be its Coxeter system.
Let $\preceq$ be the relation on $W \times \SS^f$ defined as follows. 
\[
(u,X) \preceq (v,Y)
\]
if
\[
X \subset Y\,,\ v^{-1}u \in W_Y\,, \text{ and } v^{-1}u \text{ is } (\emptyset,X) \text{-minimal}.
\]
Then $\preceq$ is a partial order relation.}

\bigskip\noindent
{\bf Proof.}
Let $(u,X) \in W \times \SS^f$.
We have $X \subset X$, $u^{-1} u = 1 \in W_X$, and $u^{-1}u=1$ is $(\emptyset,X)$-minimal, thus $(u,X) \preceq (u,X)$.

\bigskip\noindent
Let $(u,X), (v,Y) \in W \times \SS^f$ such that $(u,X) \preceq (v,Y)$ and $(v,Y) \preceq (u,X)$.
We have $X \subset Y \subset X$, thus $X=Y$.
We have $v^{-1}u \in W_X$ and $v^{-1}u$ is $(\emptyset,X)$-minimal.
But, the only $(\emptyset,X)$-minimal element lying in $W_X$ is $1$, thus $v^{-1}u=1$, therefore $v=u$.

\bigskip\noindent
Let $(u,X), (v,Y), (w,Z) \in W \times \SS^f$ such that $(u,X) \preceq (v,Y) \preceq (w,Z)$.
We have $X \subset Y$ and $Y \subset Z$, thus $X \subset Z$.
We also have $w^{-1}v \in W_Z$ and $v^{-1}u \in W_Y \subset W_Z$, thus $w^{-1}u = w^{-1}v\, v^{-1}u \in W_Z$.
Now, because $w^{-1}v$ is $(\emptyset,Y)$-minimal and $v^{-1}u$ is $(\emptyset,X)$-minimal, by Proposition 2.9, for all $u_0 \in W_X$ we have the following equalities. 
\begin{gather*}
\lg(w^{-1}u\,u_0) = \lg(w^{-1}v\, v^{-1}u\, u_0) = \lg(w^{-1}v) +\lg(v^{-1}u\, u_0) = \lg(w^{-1}v) +\lg(v^{-1}u) +\lg(u_0)\\
 = \lg(w^{-1}v\, v^{-1}u) +\lg(u_0) = \lg(w^{-1}u) +\lg(u_0)\,.
\end{gather*}
By Proposition 2.9 it follows that $w^{-1}u$ is $(\emptyset,X)$-minimal.
So, $(u,X) \preceq (w,Z)$.
\qed

\bigskip\noindent
Let $\Gamma$ be a Coxeter graph, and let $(W,S)$ be its Coxeter system.
The {\it Salvetti complex} of $\Gamma$, denoted by $\Sal(\Gamma)$, is defined to be the geometric realization of the derived complex of $(W \times \SS^f, \preceq)$. 
Note that the action of $W$ on $W \times \SS^f$ defined by $w \cdot (u,X) = (wu,X)$, $w \in W$ and $(u,X) \in W \times \SS^f$, preserves the ordering.
Hence, it induces an action of $W$ on $\Sal (\Gamma)$.

\bigskip\noindent
Now, we take a Vinberg system $(W,S)$ and we denote by $\Gamma$ the Coxeter graph of $(W,S)$, viewed as a Coxeter system. 
We go back to the notations and definitions of Subsection 2.3.
So, $\AA$ denotes the Coxeter arrangement of $(W,S)$, $C_0$ denotes the fundamental chamber of $(W,S)$, $\FF(C_0)$ denotes the set of facets $F$ of $\AA$ such that $F \preceq C_0$, and $\PP^f = \{wW_X \mid w \in W \text{ and } X \in \SS^f\}$.
Recall also that we have bijective maps $\iota : \SS^f \to \PP(C_0)$ and $\tilde \iota : \PP^f \to \FF(\AA)$ whose properties are stated in Theorem 2.10.
We define the map $\varphi: W \times \SS^f \to \Sal_0(\AA)$ as follows
\[
\varphi (w,X) = (\tilde \iota(wW_X), \tilde\iota(w)) = (w(\iota(X)), w(C_0))\,.
\]
The main result of this subsection is the following.

\bigskip\noindent
{\bf Theorem 3.3.}
{\it The map $\varphi$ is a bijective map which satisfies the following property. 
\begin{itemize}
\item[($*$)]
Let $(u,X), (v,Y) \in W \times \SS^f$.
We have $(u,X) \preceq (v,Y)$ if and only if $\varphi(v,Y) \preceq \varphi(u,X)$.
\end{itemize}}

\bigskip\noindent
Before proving Theorem 3.3, we first give two important consequences.
By construction, the map $\varphi$ induces a homeomorphism $h: \Sal(\Gamma) \to \Sal(\AA)$.
Moreover, $\varphi$ being equivariant, the homeomorphism $h$ is also equivariant.
Combining this with Theorem 3.1 we obtain the following.

\bigskip\noindent
{\bf Corollary 3.4.}
{\it There exists a homotopy equivalence $f : \Sal(\Gamma) \to M(W,S)$ equivariant under the actions of $W$ and that induces a homotopy equivalence $\bar f : \Sal(\Gamma)/W \to M(W,S)/W=N(W,S)$.}

\bigskip\noindent
The following result is a direct consequence of Corollary 3.4.
It was previously proved by Charney and Davis \cite{ChaDav1}.

\bigskip\noindent
{\bf Corollary 3.5.}
{\it The homotopy type of $N(W,S)$ (resp. $M(W,S)$) depends only on the Coxeter graph $\Gamma$.}

\bigskip\noindent
{\bf Proof of Theorem 3.3.}
Let $(F,C) \in \Sal_0(\AA)$.
There exists $w \in W$ such that $C=w(C_0)$.
We have $w^{-1}(F) \preceq w^{-1}(C)=C_0$, thus $w^{-1}(F) \in \FF(C_0)$, therefore there exists $X \in \SS^f$ such that $w^{-1}(F) = \iota(X)$.
Hence, $(F,C) = \varphi(w,X)$.
This shows that $\varphi$ is a surjective map.

\bigskip\noindent
Let $(u,X), (v,Y) \in W \times \SS^f$ such that $\varphi(u,X)=\varphi(v,Y)=(F,C)$.
We have $u(C_0) = v(C_0)=C$, thus $u=v$.
Moreover, $\iota(X) = \iota(Y) = u^{-1}(F)$, thus $X=Y$.
This shows that $\varphi$ is injective.

\bigskip\noindent
Let $(u,X), (v,Y) \in W \times \SS^f$ such that $(u,X) \preceq (v,Y)$.
Set $(F,C) = \varphi(u,X)$ and $(G,D) = \varphi(v,Y)$.
We have $X \subset Y$ and $u$ can be written in the form $u=vw$, where $w \in W_Y$ and $w$ is $(\emptyset,X)$-minimal.
We have $\iota(Y) \preceq \iota(X)$ by Theorem 2.10.I, thus $\iota(Y) = w^{-1} (\iota(Y)) \preceq w^{-1} (\iota(X))$, therefore $v(\iota(Y)) = G \preceq u(\iota(X)) = F$.
Let $H \in \AA_F$.
We have $u^{-1}(H) \in \AA_{\iota(X)}$.
By Theorem~2.10.I, there exists a reflection $r$ lying in $W_X$ such that $u^{-1}(H)=H_r$.
By Theorem~2.10.III, $u^{-1}(H)$ does not separate $C_0$ and $w^{-1}(C_0)$, thus $H$ does not separate $u(C_0)=C$ and $uw^{-1}(C_0) = v(C_0) = D$.
It follows that $C_F=D_F$, thus $D_G \subset C_F$.
Hence, $(G,D) \preceq (F,C)$.

\bigskip\noindent
Let $(u,X), (v,Y) \in W \times \SS^f$.
We set $(F,C) = \varphi (u,X)$ and $(G,D) = \varphi(v,Y)$, and we assume that $(G,D) \preceq (F,C)$.
Let $w =v^{-1}u$ (thus $u=vw$).
Since $G = \tilde \iota(vW_Y) \preceq F = \tilde \iota (uW_X)$, by Theorem 2.10.II, $uW_X \subset vW_Y$, thus, again by Theorem 2.10.II, $X \subset Y$ and $u \in vW_Y$.
The later inclusion implies that $w \in W_Y$.
Let $r$ be a reflection lying in $W_X$.
We have $H_r \in \AA_{\iota(X)}$ by Theorem 2.10.I, thus $u(H_r) \in \AA_F$.
Since $C_F = D_F$, $u(H_r)$ does not separate $C$ and $D$, thus $H_r$ does not separate $u^{-1}(C) = C_0$ and $u^{-1}(D) = w^{-1}(C_0)$.
By Theorem 2.10.III, it follows that $w$ is $(\emptyset,X)$-minimal.
We conclude that $(u,X) \preceq (v,Y)$.
\qed

\subsection{Cellular decompositions and fundamental groups}

\bigskip\noindent
Let $\Gamma$ be a Coxeter graph, and let $(W,S)$ be its Coxeter system.
For all $s \in S$ we set $W^s = W_{S \setminus \{s\}}$.
The {\it Coxeter complex} of $\Gamma$, denoted by $\Cox = \Cox (\Gamma)$, is the simplicial complex defined as follows.
\begin{itemize}
\item[(a)]
The set of vertices of $\Cox$ is the set of cosets $\{wW^s \mid w \in W \text { and } s \in S\}$.
\item[(b)]
A family $\{w_0W^{s_0}, w_1 W^{s_1}, \dots, w_p W^{s_p} \}$ is a simplex of $\Cox$ if the intersection $w_0 W^{s_0} \cap \cdots \cap w_p W^{s^p}$ is nonempty.
\end{itemize}
For $X \subset S$, $X \neq \emptyset$, we set $W^X=W_{S \setminus X}$.
Let $X \subset S$, $X \neq \emptyset$, and let $w \in W$.
With the coset $wW^X$ we associate the simplex $\Delta(wW^X) = \{wW^s \mid s \in X\}$ of $\Cox$.
Every simplex has this form, and we have $\Delta(uW^X) = \Delta(vW^Y)$ if and only if $uW^X=vW^Y$.
The Coxeter complexes play a prominent role in the definition and the study of Tits buildings.
We refer to \cite{AbrBro1} for a detailed study of these complexes and their applications to Tits buildings.
In this paper we will use the following result.
This is well-known and can be found for instance in \cite{AbrBro1}.

\bigskip\noindent
{\bf Proposition 3.6.}
{\it Let $\Gamma$ be a Coxeter graph of spherical type.
Then $|Cox(\Gamma)|$ is (homeomorphic to) a sphere of dimension $|S|-1$.}

\bigskip\noindent
Set $\PP_0=\PP_0(\Gamma) = \{wW_X \mid w \in W,\ X \subset S \text{ and } X \neq S\} = \{w W^Y \mid w \in W,\ Y \subset S \text{ and } Y \neq \emptyset\}$, and $\PP=\PP(\Gamma) = \{wW_X \mid w \in W \text{ and } X \subset S\}$ that we order by the inclusion.
Note that, if $\Gamma$ is of spherical type, then $\PP$ coincides with the set $\PP^f$ defined in Subsection 2.3.
Observe also that $|\PP'|$ is the cone over $|\PP_0'|$.
On the other hand, by the above, $\PP_0$ is (isomorphic to) the barycentric subdivision of $|\Cox(\Gamma)|$.
Then the following result follows from Proposition 3.6.

\bigskip\noindent
{\bf Corollary 3.7.}
{\it Let $\Gamma$ be a Coxeter graph of spherical type.
Then the geometric realization $|\PP'|$ of the derived complex of $\PP$ is homeomorphic to a disk of dimension $|S|$, whose boundary is the geometric realization $|\PP_0'|$ of the derived complex of $\PP_0$.}

\bigskip\noindent
There is a ``geometric'' way to describe the Coxeter complex and see Proposition 3.6. 
Recall the construction of the canonical representation (see Subsection 2.2).
We take an abstract set $\{e_s \mid s \in S\}$ in one-to-one correspondence with $S$, and we denote by $V$ the real vector space with basis $\{e_s \mid s \in S\}$.
There is a symmetric bilinear form $B: V \times V \to \R$, and a faithful linear representation $\rho : W \to \GL(V)$ that leaves invariant the form $B$ and which is called canonical representation.

\bigskip\noindent
Assume that $\Gamma$ is of spherical type.
Then $B$ is positive definite (see Theorem 2.6), thus we can identify $V^*$ with $V$ via the form $B$.
For all $s \in S$, we set $H_s = \{x \in V \mid B(x,e_s)=0\}$.
This is the hyperplane orthogonal to $e_s$.
Then $\rho(s) = \rho^*(s)$ is the orthogonal reflection with respect to $H_s$ for all $s \in S$, and $W$, identified with $\rho(W)$, is a (finite) linear group generated by reflections.
In this case the Tits cone is $I=V$.
Recall that $\AA$ denotes the set of reflection hyperplanes of $W$.
This is a finite hyperplane arrangement in $V$.
We observe that the arrangement $\AA$ defines a cellular decomposition of the sphere $\S^{|S|-1} = \{ x \in V \mid B(x,x)=1 \}$.
The proof of the following can be found in \cite{AbrBro1}.

\bigskip\noindent
{\bf Proposition 3.8.}
{\it The cellular decomposition of $\S^{|S|-1}$ determined by $\AA$ is a simplicial decomposition which is isomorphic to $\Cox(\Gamma)$.}

\bigskip\noindent
{\bf Example.}
Let $m \in \N$, $m \ge 2$.
We identify $\R^2$ with $\C$, we denote by $H_0$ the (real) line spanned by $1$, by $H_1$ the line spanned by $e^{i \pi/m}$, by $s$ the orthogonal reflection with respect to $H_0$, and by $t$ the orthogonal reflection with respect to $H_1$.
Let $W$ be the group generated by $s$ and $t$.
Then $W$ is the dihedral group of order $2m$ and has the presentation
\[
W= \langle s,t \mid s^2 = t^2 = (st)^m=1 \rangle\,.
\]
For $k \in \{0,1, \dots, m-1\}$, we denote by $H_k$ the line spanned by $e^{ik \pi/m}$.
Then $\AA=\{H_k \mid 0 \le k \le m-1\}$.
The cellular decomposition of $\S^1$ defined by $\AA$ is composed by $2m$ vertices and $2m$ edges (see Figure 3.3).
Let $\bar C_0$ be the closed cone spanned by $1$ and $e^{i \pi/m}$.
Then $a_0=\bar C_0 \cap \S^1$ is an edge of the decomposition, and $W$ acts freely and transitively on the set of edges.
For $w \in W$, the $1$-simplex of $\Cox$ corresponding to $w \, a_0$ is $\Delta(wW^{\{s,t\}})= \Delta(w\cdot \{1\})$.
The vertices adjacent to the edge $w\, a_0$ correspond to the vertices $wW^{s}=w\cdot \{1,t\}$ and $wW^{t} = w\cdot\{1,s\}$ of $\Cox$.

\begin{figure}[tbh]
\bigskip
\centerline{
\setlength{\unitlength}{0.5cm}
\begin{picture}(12,10)
\put(3,1){\includegraphics[width=4cm]{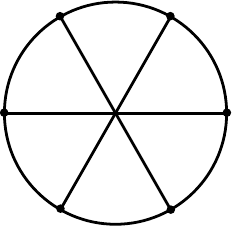}}
\put(1,1.7){\small $= tst\, a_0$}
\put(1.7,2.5){\small $sts\,a_0$}
\put(1.9,7){\small $ts\,a_0$}
\put(6.1,0.4){\small $st\,a_0$}
\put(6.5,9.2){\small $t\,a_0$}
\put(7,7){\small $H_1$}
\put(8.5,5.3){\small $H_0$}
\put(10.4,2.7){\small $s\,a_0$}
\put(10.5,6.7){\small $a_0$}
\end{picture}} 

\bigskip
\centerline{{\bf Figure 3.3.} Coxeter complex.}
\end{figure}

\bigskip\noindent
We turn back to the assumption that $\Gamma$ is any Coxeter graph. 
Let $(u,X) \in W \times \SS^f$.
Set
\[
C(u,X) = \{ (v,Y) \in W \times \SS^f \mid (v,Y) \preceq (u,X) \}\,.
\]
Furthermore, let $wW_Y \in \PP(\Gamma_X)$, and let $w_0$ be the $(\emptyset,Y)$-minimal element lying in $wW_Y$.
With the coset $wW_Y$ we associate the element $(uw_0,Y)$ of $W \times \SS^f$, that we denote by $f(wW_Y)$. 
The key point in the construction of the cellular decomposition of $\Sal (\Gamma)$ is the following.

\bigskip\noindent
{\bf Lemma 3.9.}
{\it Let $(u,X) \in W \times \SS^f$. 
Then $f(wW_Y) \in C(u,X)$ for all $wW_Y \in \PP(\Gamma_X)$, and the map $f: \PP(\Gamma_X) \to C(u,X)$ is a poset isomorphism.}

\bigskip\noindent
{\bf Proof.}
Let $wW_Y \in \PP(\Gamma_X)$.
We can assume without loss of generality that $w$ is $(\emptyset,Y)$-minimal.
Since $Y \subset X$, $w \in W_X$, and $w$ is $(\emptyset,Y)$-minimal, we have $f(wW_Y) = (uw,Y) \preceq (u,X)$.
So, $f(wW_Y) \in C(u,X)$.

\bigskip\noindent
If $(v,Y) \in C(u,X)$, then $(v,Y) = f(u^{-1}v W_Y)$.
So, $f$ is surjective.
On the other hand, if $f(wW_Y) = (v',Y')$, then $Y=Y'$ and $wW_Y = u^{-1}v'W_Y$, thus $wW_Y$ is entirely determined by its image.
So, $f$ is injective.

\bigskip\noindent
Let $w_1W_{Y_1}, w_2W_{Y_2} \in \PP(\Gamma_X)$.
We can assume without loss of generality that $w_1$ is $(\emptyset,Y_1)$-minimal and $w_2$ is $(\emptyset,Y_2)$-minimal.
Suppose $w_1 W_{Y_1} \subset w_2 W_{Y_2}$.
Then $Y_1 \subset Y_2$ and $w_1 \in w_2 W_{Y_2}$ (see Theorem 2.10.II).
By Proposition 2.9, $w_1$ is of the form $w_1 = w_2 w_1'$, where $\lg(w_1) = \lg(w_2) + \lg(w_1')$, and $w_1' \in W_{Y_2}$.
Moreover, $w_1'$ is $(\emptyset,Y_1)$-minimal because $w_1$ is.
Since $f(w_1W_{Y_1}) = (uw_1,Y_1) = (uw_2w_1',Y_1)$ and $f(w_2W_{Y_2}) = (uw_2,Y_2)$, it follows that $f(w_1 W_{Y_1}) \preceq f(w_2 W_{Y_2})$.
Suppose that $f(w_1 W_{Y_1}) \preceq f(w_2 W_{Y_2})$.
Then $Y_1 \subset Y_2$ and $w_2^{-1} w_1 \in W_{Y_2}$, thus $w_1 W_{Y_1} \subset w_2 W_{Y_2}$.
\qed

\bigskip\noindent
We describe the cellular decomposition of $\Sal (\Gamma)$ as follows.
For all $w \in W$ we have a vertex $x(w)$ corresponding to the poset $w\, \{1\}$.
The $0$-skeleton of $\Sal(\Gamma)$ is $\{x(w) \mid w \in W\}$.
For $p \in \N$, the set of $p$-cells of $\Sal (\Gamma)$ is $\{ |C(u,X)'| \mid (u,X) \in W \times \SS^f \text{ and } |X|=p\}$, and the $p$-skeleton $\Sal(\Gamma)_p$ is the union of these cells.
Lemma 3.9 and Corollary 3.7 imply that this defines a regular cellular decomposition of $\Sal(\Gamma)$.

\bigskip\noindent
For $(u,X) \in W \times \SS^f$ with $|X|=p$, we denote by $\B(u,X) = |C(u,X)'|$ the cell of $\Sal(\Gamma)_p$ associated to $(u,X)$, and we denote by $\varphi_{u,X} : \partial \B(u,X) \to \Sal(\Gamma)_{p-1}$ the gluing map. 

\bigskip\noindent
We denote by $\BSal(\Gamma)$ the quotient of $\Sal(\Gamma)$ under the action of $W$. 
Then the cellular decomposition of $\Sal(\Gamma)$ determines a cellular decomposition of $\BSal (\Gamma)$ that is described as follows.
Let $X \in \SS^f$.
The orbit of the cell $\B(1,X)$ under the action of $W$ is $\{\B(u,X) \mid u \in W\}$.
With this orbit we associate a cell $\bar \B(X)$ of $\BSal(\Gamma)$ of dimension $|X|$ and homeomorphic to $\B(1,X)$ via a homeomorphism $h_X: \bar\B(X) \to \B(1,X)$.
The set of cells of $\BSal(\Gamma)$ of dimension $p$ is $\{\bar \B(X) \mid X \in \SS^f \text{ and } |X|=p\}$.
For $X \in \SS^f$ such that $|X|=p$, the gluing map $\bar \varphi_X : \partial \bar \B(X) \to \BSal(\Gamma)_{p-1}$ is defined as follows.
\[
\bar \varphi_X = \pi \circ \varphi_{1,X} \circ h_X : \partial \bar \B(X) \to \BSal(\Gamma)_{p-1}\,, 
\]
where $\pi: \Sal(\Gamma) \to \BSal(\Gamma)$ denotes the natural projection.
Note that $\bar \varphi_X$ is not in general a homeomorphism onto its image, thus $\BSal(\Gamma)$ is not a regular CW-complex.
Note also that $\B(1,X)$ can be viewed as embedded into $\Sal(\Gamma)$, but $\bar \B (X)$ cannot be viewed as embedded into $\BSal(\Gamma)$.

\bigskip\noindent
For practical reasons (in particular, for calculating fundamental groups), and in order to better understand these complexes, we turn now to describe the $p$-skeletons of $\Sal(\Gamma)$ and $\BSal(\Gamma)$ for $p=0,1,2$.

\bigskip\noindent
{\it $0$-skeleton.}
As mentioned before, the $0$-skeleton of $\Sal(\Gamma)$ is a set $\{x(w) \mid w \in W\}$ in one-to-one correspondence with $W$.
The $0$-skeleton of $\BSal(\Gamma)$ is reduced to a point that we denote by $x_0$.

\bigskip\noindent
{\it $1$-skeleton.}
With every $(u,s) \in W \times S$ is associated an edge $\B(u,\{s\})$ of $\Sal(\Gamma)$ whose extremities are $x(u)$ and $x(us)$.
We denote this edge by $a(u,s)$, and we assume it to be oriented from $x(u)$ to $x(us)$.
So, for $u,v \in W$, if $v$ is of the form $v=us$ with $s \in S$, there is an edge $a(u,s)$ going from $x(u)$ to $x(v)$, and there is another edge $a(v,s)$ going from $x(v)$ to $x(u)$  (see Figure~3.4).
On the other hand, there is no edge joining $x(u)$ and $x(v)$ if $v$ is not of the form $v=us$ with $s \in S$. 
With every $s \in S$ is associated an edge $\bar a_s =\bar \B(\{s\})$ of $\BSal(\Gamma)$ whose extremities are both equal to $x_0$.
Let $s \in S$.
It is easily seen that the action of $W$ on $\{a(u,s) \mid u \in W\}$ preserves the orientations of the $a(u,s)$, thus it induces an orientation on $\bar a_s$.
So, we assume $\bar a_s$ to be endowed with this orientation.

\begin{figure}[tbh]
\bigskip
\centerline{
\setlength{\unitlength}{0.5cm}
\begin{picture}(17,4)
\put(2,1){\includegraphics[width=7cm]{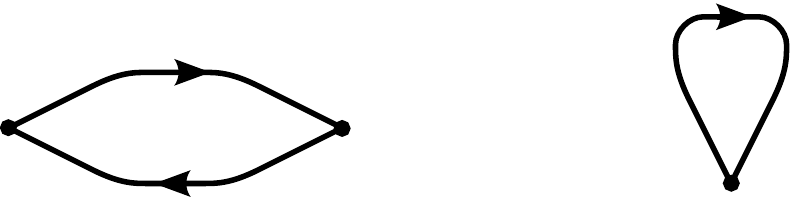}}
\put(0.5,2.2){\small $x(u)$}
\put(4,0.4){\small $a(v,s)$}
\put(4,3.7){\small $a(u,s)$}
\put(7.7,1.4){\small $=x(us)$}
\put(8.5,2.2){\small $x(v)$}
\put(14.6,0.6){\small $x_0$}
\put(15.9,2.2){\small $\bar a_s$}
\end{picture}} 

\bigskip
\centerline{{\bf Figure 3.4.} Edges in $\Sal(\Gamma)$ and in $\BSal(\Gamma)$.}
\end{figure}

\bigskip\noindent
{\it $2$-skeleton.}
Let $s,t \in S$, $s \neq t$.
Note that we have $\{s,t\} \in \SS^f$ if and only if $m_{s,t} \neq \infty$.
Assume $m=m_{s,t} \neq \infty$.
With every $u \in W$ is associated a $2$-cell of $\Sal(\Gamma)$, $\B (u, \{s,t\})$, whose boundary is 
\[
a(u,s) \, a(us,t) \cdots a(u\,\Pi(s,t:m-1),t)\, a(u\,\Pi(t,s:m-1),s)^{-1} \cdots a(ut,s)^{-1}\, a(u,t)^{-1}
\]
if $m$ is even, and 
\[
a(u,s) \, a(us,t) \cdots a(u\,\Pi(s,t:m-1),s)\,a(u\,\Pi(t,s:m-1),t)^{-1} \cdots a(ut,s)^{-1}\, a(u,t)^{-1}
\]
if $m$ is odd (see Figure 3.5).
The $W$-orbit of $2$-cells $\{ \B(u,\{s,t\}) \mid u \in W\}$ determines the $2$-cell $\bar \B(\{s,t\})$ of $\BSal(\Gamma)$.
By the above, the boundary curve of $\bar \B(\{s,t\})$ is 
\[
\bar a_s\, \bar a_t \cdots \bar a_t\, \bar a_s^{-1} \cdots \bar a_s^{-1}\, \bar a_t^{-1} = \Pi(\bar a_s, \bar a_t : m)\, \Pi(\bar a_t, \bar a_s :m)^{-1}
\]
if $m$ is even, and
\[
\bar a_s\, \bar a_t \cdots \bar a_s\, \bar a_t^{-1} \cdots \bar a_s^{-1}\, \bar a_t^{-1} = \Pi(\bar a_s, \bar a_t : m)\, \Pi(\bar a_t, \bar a_s :m)^{-1}
\]
if $m$ is odd (see Figure 3.5).

\begin{figure}[tbh]
\bigskip
\centerline{
\setlength{\unitlength}{0.5cm}
\begin{picture}(27,10)
\put(2,1.5){\includegraphics[width=12cm]{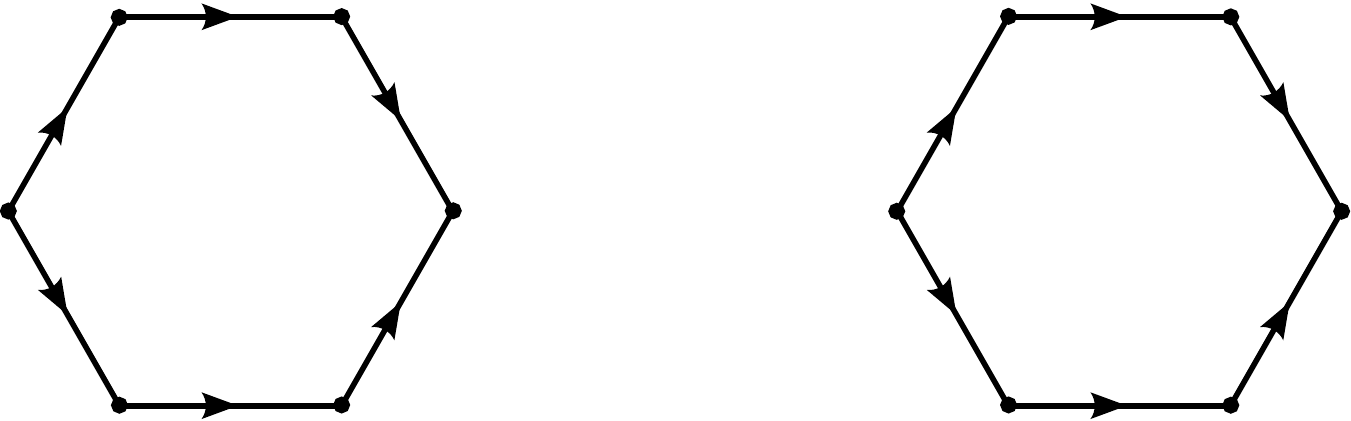}}
\put(0.5,5){\small $x(u)$}
\put(0.8,3){\small $a(u,t)$}
\put(0.8,7.1){\small $a(u,s)$}
\put(2.5,1){\small $x(ut)$}
\put(2.5,9){\small $x(us)$}
\put(4.9,0.8){\small $a(ut,s)$}
\put(4.9,9.1){\small $a(us,t)$}
\put(7.8,1){\small $x(uts)$}
\put(7.7,9){\small $x(ust)$}
\put(9.4,2.7){\small $a(uts,t)$}
\put(9.5,7){\small $a(ust,s)$}
\put(10,4.2){\small $=x(utst)$}
\put(10.5,5){\small $x(usts)$}
\put(17,5){\small $x_0$}
\put(18.3,3.2){\small $\bar a_t$}
\put(18.3,7.3){\small $\bar a_s$}
\put(19.5,1.2){\small $x_0$}
\put(19,9){\small $x_0$}
\put(21.2,9.1){\small $\bar a_t$}
\put(21.5,1){\small $\bar a_s$}
\put(23.8,1.2){\small $x_0$}
\put(24,9){\small $x_0$}
\put(25.1,2.8){\small $\bar a_t$}
\put(25.2,7){\small $\bar a_s$}
\put(26.3,5){\small $x_0$}
\end{picture}} 

\bigskip
\centerline{{\bf Figure 3.5.} $2$-cells in $\Sal(\Gamma)$ and in $\BSal(\Gamma)$.}
\end{figure}

\bigskip\noindent
A straightforward consequence of these descriptions and of Theorem 2.3 is the following.

\bigskip\noindent
{\bf Theorem 3.10.}
{\it We have $\pi_1(\BSal(\Gamma),x_0)= A_\Gamma$, $\pi_1(\Sal (\Gamma),x(1))= CA_\Gamma$, and the exact sequence associated with the regular covering $\Sal(\Gamma) \to \BSal(\Gamma)$ is the following.
\[
\xymatrix{
1 \ar[r] & CA_\Gamma \ar[r] & A_\Gamma \ar[r]^\theta & W \ar[r] & 1}\,.
\]}

\bigskip\noindent
{\bf Corollary 3.11.}
(Van der Lek \cite{Lek1}).
{\it Let $(W,S)$ be a Vinberg system.
Let $\Gamma$ be the Coxeter graph of $(W,S)$, viewed as a Coxeter system.
Then $\pi_1(N(W,S))= A_\Gamma$, $\pi_1(M(W,S))= CA_\Gamma$, and the exact sequence associated with the regular covering $M(W,S) \to N(W,S)$ is the following. 
\[
\xymatrix{
1 \ar[r] & CA_\Gamma \ar[r] & A_\Gamma \ar[r]^\theta & W \ar[r] & 1}\,.
\]}

%%%%%%%%%%%%%%
%%CHAPITRE 4%%
%%%%%%%%%%%%%%

\section{$K(\pi,1)$ problem for spherical type Artin groups}

In this section we first describe the universal cover $\TSal(\Gamma)$ of the Savetti complex of any Coxeter graph $\Gamma$ (see Subsection 4.1).
Afterwards we prove that a certain subcomplex $\TSal^+(\Gamma)$ of $\TSal(\Gamma)$ is contractible (see Subsection 4.2).
At the end, we prove that $\TSal(\Gamma)$ is contractible if $\Gamma$ is of spherical type (see Subsection 4.3).

\subsection{Universal cover of the Salvetti complex}

\bigskip\noindent
We take a Coxeter graph $\Gamma$, and we denote by $(A,\Sigma)$ the Artin system of $\Gamma$.
Recall the homomorphism $\theta : A \to W$ which sends $\sigma_s$ to $s$ for all $s \in S$.
Recall also that $\theta$ has a natural set-section $\tau: W \to A$ defined as follows (see Subsection 2.2).
Let $w \in W$.
We choose a reduced expression $\mu=s_1 \cdots s_\ell$ of $w$ and we set $\tau(w) = \sigma_{s_1} \cdots \sigma_{s_\ell}$.

\bigskip\noindent
{\bf Lemma 4.1.}
{\it Let $\preceq$ be the relation on $A \times \SS^f$ defined by  
\[
(\alpha,X) \preceq (\beta,Y)
\]
if $X \subset Y$ and $\alpha$ can be written in the form $\alpha=\beta\, \tau(w)$, where $w \in W_Y$ and $w$ is $(\emptyset,X)$-minimal.
Then $\preceq$ is a partial order relation.}

\bigskip\noindent
{\bf Proof.}
We cleary have $(\alpha,X) \preceq (\alpha,X)$ for all $(\alpha,X) \in A \times \SS^f$.
Let $(\alpha,X), (\beta, Y) \in A \times \SS^f$ such that $(\alpha,X) \preceq (\beta,Y)$ and $(\beta, Y) \preceq (\alpha,X)$.
We have $X \subset Y$ and $Y \subset X$, thus $X=Y$.
By definition, $\alpha$ can be written in the form $\alpha = \beta \, \tau(u)$, where $u \in W_X$ and $u$ is $(\emptyset,X)$-minimal.
But, the only $(\emptyset,X)$-minimal element lying in $W_X$ is $1$, thus $u=1$ and $\alpha=\beta$.

\bigskip\noindent
Let $(\alpha, X), (\beta, Y), (\gamma,Z) \in A \times \SS^f$ such that $(\alpha,X) \preceq (\beta,Y)$ and $(\beta,Y) \preceq (\gamma,Z)$.
We have $X \subset Y \subset Z$.
Moreover, $\alpha$ can be written in the form $\alpha = \beta \, \tau(u)$, where $u \in W_Y$ and $u$ is $(\emptyset,X)$-minimal, and $\beta$ can be written in the form $\beta = \gamma\, \tau(v)$, where $v \in W_Z$, and $v$ is $(\emptyset,Y)$-minimal.
Set $w=vu$.
Since $u \in W_Y$ and $v$ is $(\emptyset,Y)$-minimal, by Proposition 2.9, we have $\lg(w) = \lg(v) + \lg(u)$, thus $\tau(w) = \tau(v)\, \tau(u)$.
Hence, $\alpha = \gamma \, \tau(w)$.
We have $w \in W_Z$, since $v,u \in W_Z$.
On the other hand, one can easily prove following the same arguments as in the proof of Lemma 3.2 that $w$ is $(\emptyset,X)$-minimal.
So, $(\alpha,X) \preceq (\gamma,Z)$.
\qed

\bigskip\noindent
We denote by $\TSal(\Gamma)$ the geometric realization of the derived complex of $(A \times \SS^f, \preceq)$.
The action of $A$ on $A \times \SS^f$ defined by
\[
\beta \, (\alpha,X) = (\beta\alpha,X)
\]
induces a free and properly discontinuous action of $A$ on $\TSal(\Gamma)$.
It is easily shown that $\TSal(\Gamma)/CA = \Sal(\Gamma)$ and $\TSal(\Gamma)/A = \BSal(\Gamma)$.
Hence, since the fundamental group of $\BSal(\Gamma)$ is $A$, we have the following. 

\bigskip\noindent
{\bf Proposition 4.2.}
{\it $\TSal(\Gamma)$ is the universal cover of $\Sal(\Gamma)$ and of $\BSal(\Gamma)$.}

\subsection{The subcomplex $\TSal^+(\Gamma)$}

\bigskip\noindent
Recall that the Artin monoid of a Coxeter graph $\Gamma$ is the monoid $A_\Gamma^+$ that admits the following monoid presentation:
\[
A_\Gamma^+ = \langle \Sigma \mid \Pi(\sigma_s,\sigma_t : m_{s,t}) = \Pi(\sigma_t,\sigma_s : m_{s,t}) \text{ for all } s,t \in S,\ s \neq t,\ m_{s,t} \neq \infty \rangle^+ \,.
\]
By \cite{Paris2}, the natural homomorphism $A_\Gamma^+ \to A_\Gamma$ is injective.

\bigskip\noindent
We define $\TSal^+(\Gamma)$ to be the geometric realization of the derived complex of $(A_\Gamma^+ \times \SS^f, \preceq)$, where $\preceq$ is the restriction to $A_\Gamma^+ \times \SS^f$ of the ordering defined in Lemma 4.1.
The main result of this subsection is the following.

\bigskip\noindent
{\bf Theorem 4.3.}
{\it The subcomplex $\TSal^+(\Gamma)$ is contractible.}

\bigskip\noindent
The remainder of the subsection is dedicated to the proof of Theorem 4.3.

\bigskip\noindent
Let $V=\R^\ell$, and let $H_1, \dots, H_k$ be hyperplanes of $V$.
For each $i \in \{1, \dots, k\}$ we take a linear form $\alpha_i : V \to \R$ such that $\Ker \, \alpha_i = H_i$.
We say that $H_1, \dots, H_k$ are {\it independent} if $\alpha_1, \dots, \alpha_k$ are linearly independent in $V^*$.

\bigskip\noindent
{\bf Lemma 4.4.}
{\it Let $I$ be a nonempty open convex cone in $V=\R^\ell$, and let $H_1, \dots, H_k$ be independent hyperplanes.
Set $L=H_1 \cap \cdots \cap H_k$, and assume that $L \cap I \neq \emptyset$.
So, we have $H_i \cap I \neq \emptyset$ for all $i \in \{1, \dots, k\}$.
For each $i \in \{1, \dots, k \}$ we take an open half-space bounded by $H_i$ that we denote by $H_i^+$.
Then $(\cup_{i=1}^k H_i^+) \cap I$ is contractible.}

\bigskip\noindent
{\bf Proof.}
We choose a basis $\{e_1, \dots, e_\ell\}$ for $V$ so that $H_i$ is defined by the equality $x_i=0$ with respect to this basis, and $H_i^+$ is defined by the inequality $x_i > 0$, for all $i \in \{1, \dots, k\}$.
Choose a point $p_0 \in L \cap I$.
Since $p_0 \in L$, it can be written in the form $p_0=\lambda_{k+1}e_{k+1} + \cdots + \lambda_\ell e_\ell$, where $\lambda_{k+1}, \dots, \lambda_\ell \in \R$.
Since $I$ is open, there is $\varepsilon >0$ such that 
\[
q_0 = \varepsilon e_1 + \cdots + \varepsilon e_k + p_0 = \varepsilon e_1 + \cdots + \varepsilon e_k + \lambda_{k+1}e_{k+1} + \cdots + \lambda_\ell e_\ell \in I\,.
\]
For all $t \in [0,1]$ we define $h_t : V \to V$ by
\[
h_t(p) = (1-t)p +t\,q_0\,, \quad p \in V\,.
\]
It is easily checked that $h_t(p) \in (\cup_{i=1}^k H_i^+) \cap I$ if $p \in (\cup_{i=1}^k H_i^+) \cap I$, that $h_0(p) = p$ for all $p \in V$, that $h_1(p)=q_0$ for all $p \in V$, and that $h_t(q_0)=q_0$ for all $t \in [0,1]$.
\qed

\bigskip\noindent
Recall that $\PP^f$ denotes the set $\{wW_X \mid w \in W \text{ and } X \in \SS^f\}$ (see Subsection 2.3). 
Note that every coset $wW_X \in \PP^f$ is uniquely represented by the pair $(u,X)$, where $u$ is the unique $(\emptyset,X)$-minimal element lying in $wW_X$.
For $s \in S$, we say that the coset $wW_X$ is {\it $s$-minimal} if $lg(su) = \lg(u)+1$ and $su$ is $(\emptyset,X)$-minimal.
We denote by $\PP^f_s$ the set of cosets $wW_X$ lying in $\PP^f$ that are $s$-minimal, and we assume $\PP^f_s$ ordered by the inclusion.
For $X_0 \in \SS^f$, $X_0 \neq \emptyset$, we set
\[
\PP_{X_0}^f = \bigcup_{s \in X_0} \PP^f_s\,.
\]

\bigskip\noindent
{\bf Lemma 4.5.}
{\it Let $s \in S$, and let $uW_X, vW_Y \in \PP^f$.
If $uW_X \in \PP^f_s$ and $vW_Y \subset uW_X$, then $vW_Y \in \PP_s^f$.}

\bigskip\noindent
{\bf Proof.}
We can assume without loss of generality that $u$ is $(\emptyset,X)$-minimal and that $v$ is $(\emptyset,Y)$-minimal.
We observe that, by Proposition 2.9, we have $uW_X \in \PP_s^f$ if and only if $\lg(suw) = \lg(u) + lg(w) +1$ for all $w \in W_X$.
Since $vW_Y \subset uW_X$, we have $Y \subset X$ and $v \in uW_X$.
We write $v=uv'$, where $v' \in W_X$.
Since $u$ is $(\emptyset,X)$-minimal, we have $\lg(v) = \lg(u) + \lg(v')$.
Moreover, since $v$ is $(\emptyset,Y)$-minimal, $v'$ is also $(\emptyset,Y)$-minimal.
Let $w \in W_Y$.
Then
\[
\lg(svw) = \lg(suv'w) = \lg(u) + \lg(v'w)+1 = \lg(u) + \lg(v') + \lg(w) +1 = \lg(v)+\lg(w)+1\,.
\]
This shows that $vW_Y \in \PP_s^f$.
\qed

\bigskip\noindent
{\bf Lemma 4.6.}
{\it \begin{itemize}
\item[(1)]
The geometric realization $| (\PP^f)'|$ of the derived complex of $\PP^f$ is contractible.
\item[(2)]
Let $X_0 \in \SS^f$, $X_0 \neq \emptyset$. 
Then the geometric realization $|(\PP^f_{X_0})'|$ of the derived complex of $\PP^f_{X_0}$ is contractible.
\end{itemize}}

\bigskip\noindent
{\bf Proof.}
Let $\rho^* : W \to \GL(V^*)$ be the dual representation of the canonical representation.
Recall that, by Theorem 2.5, this representation is faithful and, $(W,S)$, identified with $(\rho^*(W), \rho^*(S))$, is a Vinberg system. 
We denote by $I$ the Tits cone, and by $\AA$ the Coxeter arrangement in $I$ associated to $(W,S)$.
For a reflection $r$ lying in $W$ we denote by $H_r$ the fix hyperplane of $r$.
Note that, by construction, the set $\{H_s \mid s \in S \}$ is independent.

\bigskip\noindent
We denote by $C_0$ the fundamental chamber, and by $\FF(C_0)$ the set of facets $F$ of $\AA$ such that $F \preceq C_0$.
Recall that we have a bijection $\iota : \SS^f \to \FF(C_0)$, and this bijection extends to a bijection $\tilde \iota:\PP^f \to \FF(\AA)$, $uW_X \mapsto u\, \iota(X)$.
Moreover, we have $uW_X \subset vW_Y$ if and only if $v\,\iota(Y) \preceq u\, \iota(X)$ (see Theorem 2.10).

\bigskip\noindent
For all $s \in S$, we denote by $H_s^+$ the open half-space bounded by $H_s$ and containing $C_0$.
For $X_0 \in \SS^f$, $X_0 \neq \emptyset$, we set
\[
\Omega(X_0) = \left( \bigcup_{s \in X_0} H_s^+ \right) \cap I\,.
\]
Note that $\iota(X_0)$ is contained in $\cap_{s \in X_0} H_s$, thus $(\cap_{s \in X_0} H_s)\cap I \neq \emptyset$, therefore, by Lemma 4.4, $\Omega(X_0)$ is contractible. 

\bigskip\noindent
{\bf Claim 1.}
{\it Let $s \in S$, and let $uW_X \in \PP^f$.
We have $uW_X \in \PP_s^f$ if and only if $u\,\iota(X)$ is contained in $H_s^+$.}

\bigskip\noindent
{\it Proof of Claim 1.}
We can assume without loss of generality that $u$ is $(\emptyset,X)$-minimal.
Set $F=u\,\iota(X)$.
Observe that $F$ is contained in $H_s^+$ if and only if $C$ is contained in $H_s^+$ for every chamber $C \in \CC(\AA)$ satisfying $F \preceq C$.

\bigskip\noindent
Suppose that $uW_X \in \PP^f_s$.
Let $C \in \CC(\AA)$ such that $F \preceq C$.
By Theorem 2.10.I, there exists $w \in W_X$ such that $C=uw(C_0)$.
Then
\[
\begin{array}{ccl}
& \lg(suw)=\lg(uw)+1 & (\text{since } uW_X\in \PP_s^f)\\
\Rightarrow & \lg(w^{-1} u^{-1} s) = \lg(w^{-1} u^{-1})+1\\
\Rightarrow & C=uw(C_0) \subset H_s^+ & (\text{by Theorem 2.10.III})
\end{array}
\]
This shows that $F$ is included in $H_s^+$.

\bigskip\noindent
Suppose now that $F$ is included in $H_s^+$.
Let $w \in W_X$.
Set $C=uw(C_0)$.
Since $F \preceq C$, we have $C \subset H_s^+$.
By Theorem 2.10.III, it follows that
\[
\lg(w^{-1} u^{-1} s) = \lg(w^{-1} u^{-1})+1 \ \Rightarrow\ \lg(suw)=\lg(uw)+1 = \lg(u) + \lg(w) +1\,.
\]
This shows that $uW_X \in \PP_s^f$.

\bigskip\noindent
For every facet $F$ of $\AA$, we denote by $\omega(F)$ the set defined in the proof of Theorem ~3.1 (see Figure 3.2).
For $uW_X \in \PP^f$, we set $\omega(uW_X) = \omega(u\, \iota(X))$.
Let $s \in S$, and let $uW_X \in \PP_s^f$.
Set $F=u \,\iota(X)$.
By Claim 1, we have $F \subset H_s^+$.
Moreover, if $G$ is a facet of $\AA$ such that $F \preceq G$, then $G \subset H_s^+$.
Since $\omega(uW_X) = \omega(F)$ is contained in the union of the facets $G$ of $\AA$ satisfying $F \preceq G$, it follows that $\omega(uW_X) \subset H_s^+$.
This proves the following.

\bigskip\noindent
{\bf Claim 2.}
{\it \begin{itemize}
\item[(1)]
The set $\{ \omega(uW_X) \mid uW_X \in \PP^f \}$ is a cover of $I$ by open subsets.
\item[(2)]
Let $X_0 \in \SS^f$, $X_0 \neq \emptyset$.
Then the set $\{ \omega(uW_X) \mid uW_X \in \PP_{X_0}^f \}$ is a cover of $\Omega(X_0)$ by open subsets.
\end{itemize}}

\bigskip\noindent
By Claims 5, 6, and 7 in the proof of Theorem 3.1, we have
\[
\omega(u_0 W_{X_0}) \cap \omega(u_1 W_{X_1}) \cap \cdots \cap \omega(u_p W_{X_p}) \neq \emptyset
\]
if and only if, up to permutation, we have
\[
u_p W_{X_p} \subset \cdots \subset u_1 W_{X_1} \subset u_0 W_{X_0}\,.
\]
Moreover, in that case, this intersection is contractible.
By Theorem 2.1, it follows that $| (\PP^f)'|$ has the same homotopy type as $I$, thus $| (\PP^f)'|$ is contractible.
Similarly, if $X_0 \in \SS^f$, $X_0 \neq \emptyset$, then $| (\PP_{X_0}^f)'|$ has the same homotopy type as $\Omega(X_0)$, thus $| (\PP_{X_0}^f)'|$ is contractible by Lemma~4.4.
\qed

\bigskip\noindent
We turn back to the universal cover $\TSal(\Gamma)$ of the Salvetti complex.
For $\alpha \in A_\Gamma$, we set $\tilde C(\alpha) = \{(\alpha\tau(u),X) \mid X \in \SS^f,\ u \in W \text{ and } u \text{ is } (\emptyset,X)\text{-minimal}\}$.
We restrict the ordering of $A_\Gamma \times \SS^f$ to $\tilde C(\alpha)$, and we denote by $\Phi(\alpha) \subset \TSal(\Gamma)$ the geometric realization of the derived complex of $\tilde C(\alpha)$.

\bigskip\noindent
Let $\alpha \in A_\Gamma$, and let $uW_X \in \PP^f$.
We can assume without loss of generality that $u$ is $(\emptyset,X)$-minimal.
Then we set $f_\alpha(uW_X) = (\alpha\tau(u),X) \in \tilde C(\alpha)$.

\bigskip\noindent
{\bf Lemma 4.7.}
{\it Let $\alpha \in A_\Gamma$.
Then the map $f_\alpha : \PP^f \to \tilde C(\alpha)$ is a poset isomorphism.
In particular, it induces a homeomorphism $f_\alpha : | (\PP^f)' | \to \Phi(\alpha)$.
Hence, by Lemma 4.6, $\Phi(\alpha)$ is contractible.}

\bigskip\noindent
{\bf Proof.}
It is easily seen that $f_\alpha$ is bijective.
Let $uW_X, vW_Y \in \PP^f$.
We assume without loss of generality that $u$ is $(\emptyset,X)$-minimal and $v$ is $(\emptyset,Y)$-minimal.
Suppose that $uW_X \subset vW_Y$.
Then $X \subset Y$ and $u \in vW_Y$.
Let $u' \in W_Y$ such that $u=vu'$.
Since $v$ is $(\emptyset,Y)$-minimal, we have $\lg(u) = \lg(v) + \lg(u')$, thus $\tau(u) = \tau(v)\,\tau(u')$.
Moreover, $u'$ is $(\emptyset,X)$-minimal, since $u$ is.
This implies that $(\alpha\, \tau(u),X) =( \alpha\, \tau(v)\,\tau(u'),X) \preceq (\alpha\,\tau(v),Y)$.

\bigskip\noindent
Suppose that $(\alpha\, \tau(u), X) \preceq (\alpha\, \tau(v),Y)$.
Then $X \subset Y$ and $\alpha\,\tau(u)$ can be written in the form $\alpha\,\tau(u) = \alpha\,\tau(v)\, \tau(u')$, where $u' \in W_Y$ and $u'$ is $(\emptyset,X)$-minimal.
The latter equality implies that $u=vu' \in vW_Y$, thus $uW_X \subset vW_Y$.
\qed

\bigskip\noindent
For $n \in \N$, we set
\[
\TSal^{(n)}(\Gamma) = \bigcup_{ \alpha \in A_\Gamma^+,\ \lg(\alpha) \le n} \Phi(\alpha)\,.
\]
Note that
\[
\TSal^+(\Gamma) = \bigcup_{ \alpha \in A_\Gamma^+} \Phi(\alpha) = \bigcup_{n=0}^\infty \TSal^{(n)}(\Gamma)\,.
\]

\bigskip\noindent
{\bf Lemma 4.8.}
{\it Let $n \in \N$, and let $\alpha, \beta \in A_\Gamma^+$ such that $\alpha \neq \beta$ and $\lg(\alpha) = \lg(\beta)=n+1$. 
Then $\Phi(\alpha) \cap \Phi(\beta) \subset \TSal^{(n)}(\Gamma)$.}

\bigskip\noindent
{\bf Proof.}
For $n \in \N$, we set
\[
\TSal_0^{(n)}(\Gamma) = \bigcup_{ \alpha \in A_\Gamma^+,\ \lg(\alpha) \le n} \tilde C(\alpha)\,.
\]
Since $\Phi(\alpha) = |\tilde C(\alpha)'|$ and $\Phi(\beta) = |\tilde C(\beta)' |$, we have
\[
\Phi(\alpha) \cap \Phi(\beta) = | (\tilde C(\alpha) \cap \tilde C(\beta))'|\,.
\]
Hence, in order to prove Lemma 4.8, it suffices to show that $\tilde C(\alpha) \cap \tilde C(\beta) \subset \TSal_0^{(n)}(\Gamma)$.

\bigskip\noindent
Let $(\gamma,Z) \in \tilde C(\alpha) \cap \tilde C(\beta)$.
There exist $u,v \in W$, both $(\emptyset,Z)$-minimal, such that $\gamma = \alpha\, \tau(u) = \beta\, \tau(v)$.
Since $\tau(u) \preceq_R \gamma$ and $\tau(v) \preceq_R \gamma$, the element $\tau(u) \vee_R \tau(v)$ exists, and, by Proposition~2.12, there exists $w \in W$ such that $\tau(u) \vee_R \tau(v) = \tau(w)$.

\bigskip\noindent
Let $\mu \in A_\Gamma^+$ such that $\tau(w)= \mu\, \tau(u)$.
Set $u' = \theta(\mu)$. 
Note that $w = \theta(\tau(w)) = u'\,u$. 
We have
\[
\lg(w) \le \lg(u') + \lg(u) \le \lg(\mu) + \lg(\tau(u)) = \lg(\tau(w)) = \lg(w)\,,
\]
thus $\mu=\tau(u')$ and $\lg(w) = \lg(u') + \lg(u)$.
Similarly, there exists $v' \in W$ such that $w=v'\, v$ and $\lg(w) = \lg(v') + \lg(v)$. 
Note that this implies that $\tau(w) = \tau(u')\, \tau(u) = \tau(v') \, \tau(v)$.

\bigskip\noindent
Let $\gamma_0 \in A_\Gamma^+$ such that $\gamma = \gamma_0\, \tau(w)$.
By left cancellation, we have $\alpha = \gamma_0 \tau(u')$ and $\beta = \gamma_0 \tau(v')$.
These two equalities imply that $\lg(\gamma_0) \le n+1$.
Moreover, if we had $\lg(\gamma_0) = n+1 = \lg(\alpha) = \lg(\beta)$, then we would have $\gamma_0 = \alpha= \beta$, which is not true.
Hence, $\lg(\gamma_0) \le n$.
So, in order to prove that $(\gamma,Z) \in \TSal_0^{(n)}(\Gamma)$, it suffices to show that $w$ is $(\emptyset,Z)$-minimal.

\bigskip\noindent
Suppose that $w$ is not $(\emptyset,Z)$-minimal.
Then, by Proposition 2.9, there exists $s \in Z$ such that $\lg(ws) < \lg(w)$.
Let $u=s_1 \cdots s_\ell$ be a reduced expression of $u$, and let $u' = t_1 \cdots t_k$ be a reduced expression of $u'$.
Note that $ t_1 \cdots t_ks_1 \cdots s_\ell$ is a reduced expression of $w=u'u$.
By Theorem~2.8, either there exists $i \in \{1, \dots, \ell\}$ such that $ws= t_1 \cdots t_k s_1 \cdots \widehat{s_i} \cdots s_\ell$, or there exists $j \in \{1, \dots, k\}$ such that $ws= t_1 \cdots \widehat{t_j} \cdots t_k s_1 \cdots s_\ell$.
But, we cannot have $ws= t_1 \cdots t_k s_1 \cdots \widehat{s_i} \cdots s_\ell$, because $\lg(us) = \lg(u)+1$ (since $u$ is $(\emptyset,Z)$-minimal), thus we have $ws= t_1 \cdots \widehat{t_j} \cdots t_k s_1 \cdots s_\ell$ for some $j \in \{1, \dots, k\}$.
Set $u''=t_1 \cdots \widehat{t_j} \cdots t_k$.
Then $ws = u''u$ and $\lg(ws) = \lg(u'') + \lg(u) = \lg(w)-1$.
In particular, $\tau(ws) = \tau(u'')\, \tau(u)$.
Similarly, there exists $v'' \in W$ such that $ws=v''v$ and $\lg(ws) = \lg(v'') + \lg(v)$ (hence, $\tau(ws) = \tau(v'') \, \tau(v)$).
This contradicts the fact that $\tau(w) = \tau(u) \vee_R \tau(v)$, since $\lg(\tau(ws)) < \lg(\tau(w))$, and, by the above, $\tau(u), \tau(v) \preceq_R \tau(ws)$.
\qed

\bigskip\noindent
For $\alpha \in A_\Gamma^+$, we set
\[
\End(\alpha) = \{ s \in S \mid \sigma_s \preceq_R \alpha \}\,.
\]
Note that $\alpha$ is an upper bound of $\Sigma_{\End(\alpha)}$ for the relation $\preceq_R$, thus, by Proposition 2.12, we have $\End(\alpha) \in \SS^f$.

\bigskip\noindent
{\bf Lemma 4.9.}
{\it Let $\alpha \in A_\Gamma^+$ such that $\lg(\alpha)=n+1$.
Set $X_0 = \End(\alpha)$.
Then 
\[
\Phi(\alpha) \cap \TSal^{(n)} (\Gamma) = f_\alpha (| (\PP_{X_0}^f)'|)\,.
\]
So, by Lemma 4.6, $\Phi(\alpha) \cap \TSal^{(n)} (\Gamma)$ is contractible.}

\bigskip\noindent
{\bf Proof.}
Since $\Phi (\alpha) = | \tilde C(\alpha)' |$ and $\TSal^{(n)}(\Gamma) = | (\TSal_0^{(n)}(\Gamma))'|$, we have
\[
\Phi(\alpha) \cap \TSal^{(n)}(\Gamma) = | (\tilde C(\alpha) \cap \TSal_0^{(n)}(\Gamma))'|\,.
\]
Hence, it suffices to show the following. 
\[
\tilde C(\alpha) \cap \TSal_0^{(n)} (\Gamma) = f_\alpha (\PP_{X_0}^f)\,.
\]

\bigskip\noindent
Let $uW_X \in \PP_{X_0}^f$.
As ever, we assume that $u$ is $(\emptyset,X)$-minimal.
By definition, there exists $s \in X_0$ such that $\lg(su) = \lg(u) +1$ and $su$ is $(\emptyset,X)$-minimal.
Since $X_0=\End(\alpha)$, we have $\sigma_s \preceq_R \alpha$, thus $\alpha$ can be written in the form $\alpha = \beta \sigma_s$, where $\beta \in A_\Gamma^+$.
Since $\lg(su) = \lg(u)+1$, we have $\tau(su) = \sigma_s\,\tau(u)$, thus $\alpha\, \tau(u) = \beta\,\tau(su)$, therefore $f_\alpha (uW_X) = (\alpha\, \tau(u),X) = (\beta\,\tau(su),X) \in \tilde C(\beta)$.
On the other hand, $\lg(\beta)=n$, thus $f_\alpha (uW_X) \in \tilde C(\alpha) \cap \TSal_0^{(n)}(\Gamma)$.

\bigskip\noindent
Let $(\gamma,Z) \in \tilde C(\alpha) \cap \TSal_0^{(n)}(\Gamma)$.
There exists $\beta \in A_\Gamma^+$ such that $\lg(\beta) \le n$ and $(\gamma,Z) \in \tilde C(\beta)$.
There exist $u,v \in W$ such that $u$ and $v$ are $(\emptyset,Z)$-minimal and $\gamma = \alpha\, \tau(u) = \beta\, \tau(v)$.
It is easily shown in the same way as in the proof of Lemma 4.8 that $\tau(u) \vee_R \tau(v)$ exists, that this element is of the form $\tau(w)$ with $w \in W$, that $w$ is $(\emptyset,Z)$-minimal, and that $w$ can be written in the form $w=u'u = v'v$ with $\lg(w) = \lg(u')+\lg(u) = \lg(v')+\lg(v)$.
Note that, since $\lg(\beta) < \lg(\alpha)$ and $\alpha\,\tau(u) = \beta\,\tau(v)$, we have $\lg(u)<\lg(v)$, thus $u' \neq 1$.
We choose $s \in S$ such that $\lg(u's) = \lg(u')-1$.
Set $u''= u's$.
By the above, $\tau(u') = \tau(u'')\, \sigma_s$.
Let $\gamma_0 \in A_\Gamma^+$ such that $\gamma=\gamma_0\,\tau(w)$.
We have $\alpha = \gamma_0 \, \tau(u') = \gamma_0\, \tau(u'')\, \sigma_s$, thus $s \in X_0$.
Finally, $\lg(su) = \lg(u) +1$ (since $\lg(w) = \lg(u') + \lg(u)$) and $su$ is $(\emptyset,Z)$-minimal (since $w$ is $(\emptyset,Z)$-minimal), thus $uW_Z\in\PP_{X_0}^f$.
This shows that $(\gamma,Z) = f_\alpha(uW_Z) \in f_\alpha( \PP_{X_0}^f)$.
\qed

\bigskip\noindent
{\bf Proof of Theorem 4.3.}
We start by proving by induction on $n$ that $\TSal^{(n)}(\Gamma)$ is contractible.
If $n=0$, then $\TSal^{(0)}(\Gamma)=\Phi(1)$, thus, by Lemma 4.7, $\TSal^{(0)}(\Gamma)$ is contractible.

\bigskip\noindent
Now, we assume that $n \ge 0$ and $\TSal^{(n)}(\Gamma)$ is contractible.
Let $\alpha \in A_\Gamma^+$ such that $\lg(\alpha)=n+1$.
By Lemma 4.7, $\Phi(\alpha)$ is contractible.
Moreover, by Lemma 4.9, $\Phi(\alpha) \cap \TSal^{(n)} (\Gamma)$ is contractible.
It follows that the embedding of $\Phi(\alpha) \cap \TSal^{(n)} (\Gamma)$ into $\Phi(\alpha)$ is a homotopy equivalence, thus $\Phi(\alpha) \cap \TSal^{(n)} (\Gamma)$ is a deformation retract of $\Phi(\alpha)$ (see \cite[Thm.~4.5]{Hatch1}).
We fix a deformation retraction $h_\alpha : \Phi(\alpha) \times [0,1] \to \Phi(\alpha)$ of $\Phi(\alpha)$ onto $\Phi(\alpha) \cap \TSal^{(n)} (\Gamma)$.

\bigskip\noindent
We define a map
\[
h : \TSal^{(n+1)}(\Gamma) \times [0,1] \to \TSal^{(n+1)}(\Gamma)
\]
as follows.
Let $\alpha \in A_\Gamma^+$ such that $\lg(\alpha) \le n+1$. 
If $\lg(\alpha) \le n$, we set $h(x,t) = x$ for all $(x,t) \in \Phi(\alpha) \times [0,1]$.
If $\lg(\alpha)=n+1$, we set $h(x,t)=h_\alpha(x,t)$ for all $(x,t) \in \Phi(\alpha) \times [0,1]$.
Lemma 4.8 implies that $h$ is well-defined.
It is clear from the above that $h$ is a deformation retraction of $\TSal^{(n+1)} (\Gamma)$ onto $\TSal^{(n)} (\Gamma)$, thus $\TSal^{(n+1)} (\Gamma)$ is contractible as $\TSal^{(n)} (\Gamma)$ is.

\bigskip\noindent
Since 
\[
\TSal^+ (\Gamma) = \bigcup_{n=0}^\infty \TSal^{(n)} (\Gamma)\,,
\]
we conclude by Theorem 2.1 that $\TSal^+(\Gamma)$ is contractible.
\qed

\subsection{$K(\pi,1)$ problem for Artin groups of spherical type}

\bigskip\noindent
{\bf Theorem 4.10}
(Deligne \cite{Delig1}).
{\it If $\Gamma$ is a spherical type Coxeter graph, then $\Sal(\Gamma)$ is an Eilenberg MacLane space.}

\bigskip\noindent
{\bf Proof.}
According to the statement of Theorem 2.13, we set $\Delta = \vee_L \Sigma$.
We have the following chain of subcomplexes.
\[
\TSal^+(\Gamma) \subset \Delta^{-1}\, \TSal^+(\Gamma) \subset \cdots \subset \Delta^{-n}\, \TSal^+(\Gamma) \subset \Delta^{-n-1}\, \TSal^+(\Gamma) \subset \cdots\,.
\]
The subcomplex $\Delta^{-n}\,\TSal^+(\Gamma)$ is contractible by Theorem 4.3, and $\cup_{n=0}^\infty \Delta^{-n}\, \TSal^+(\Gamma) = \TSal(\Gamma)$ by Theorem 2.13.
We conclude by Theorem 2.1 that $\TSal(\Gamma)$ is contractible.
\qed

\section{Parabolic subgroups, FC type Artin groups, and generalizations}

Let $\Gamma$ be a Coxeter graph, and let $(A_\Gamma,\Sigma)=(A,\Sigma)$ be its Artin system.
For $X \subset S$, we set $\Sigma_X=\{\sigma_s \mid s \in X\}$, and we denote by $A_X$ the subgroup of $A$ generated by $\Sigma_X$.
Such a subgroup is called {\it standard parabolic subgroup} of $A$.
Recall that, for $X\subset S$, we set $M_X = (m_{s,t})_{s,t \in X}$, where $M=(m_{s,t})_{s,t \in S}$ is the Coxeter matrix of the Coxeter graph $\Gamma$, we denote by $\Gamma_X$ the Coxeter graph of $M_X$, and we denote by $W_X$ the subgroup of $W=W_\Gamma$ generated by $X$.
By \cite{Bourb1}, the pair $(W_X,X)$ is the Coxeter system of $\Gamma_X$.
The subgroup $W_X$ is called {\it standard parabolic subgroup} of $W$.

\bigskip\noindent
Let $T$ be a subset of $S$. 
Set $\SS^f_T=\{X \in \SS^f \mid X \subset T\}$.
Observe that the inclusion $(W_T \times \SS^f_T) \hookrightarrow (W \times \SS^f)$ induces an embedding $\iota_T:\Sal (\Gamma_T) \hookrightarrow \Sal(\Gamma)$ which is equivariant under the action of $W_T$. 
The starting point of the present section is the following theorem proved in \cite{GodPar1}.
It will be the key in the proofs of several results on standard parabolic subgroups and on some families of Artin groups. 

\bigskip\noindent
{\bf Theorem 5.1}
(Godelle, Paris \cite{GodPar1}).
{\it Let $T$ be a subset of $S$. 
Then the embedding $\iota_T:\Sal (\Gamma_T) \hookrightarrow \Sal (\Gamma)$ admits a retraction $\pi_T: \Sal (\Gamma) \to \Sal (\Gamma_T)$ which is equivariant under the action of $W_T$.}

\bigskip\noindent
{\bf Proof.}
It suffices to determine a function $\pi_T : (W \times \SS^f) \to (W_T \times \SS^f_T)$ that satisfies the following properties.
\begin{itemize}
\item
$\pi_T(u,X)= (u,X)$ for all $(u,X) \in W_T \times \SS^f_T$,
\item
$\pi_T$ is equivariant under (left) action of $W_T$,
\item
if $(u,X) \preceq (v,Y)$, then $\pi_T(u,X) \preceq \pi_T(v,Y)$.
\end{itemize}

\bigskip\noindent
Let $(u,X) \in W \times \SS^f$. 
We write $u=u_0u_1$, where $u_0 \in W_T$ and $u_1$ is $(T,\emptyset)$-minimal.
Let $X_0 = T \cap u_1 X u_1^{-1}$. 
Then we set  
\[
\pi_T(u,X) = (u_0,X_0)\,.
\]
Note that, since $W_{X_0} \subset u_1 W_X u_1^{-1}$, the group $W_{X_0}$ is finite, thus $X_0 \in \SS_T^f$.

\bigskip\noindent
It is easily seen that $\pi_T(u,X)=(u,X)$ for all $(u,X) \in W_T \times \SS^f_T$, and that $\pi_T$ is equivariant under the action of $W_T$. 
So, it remains to show that, if $(u,X), (v,Y) \in W \times \SS^f$ are such that $(u,X) \preceq (v,Y)$, then $\pi_T(u,X) \preceq \pi_T(v,Y)$.

\bigskip\noindent
Let $(u,X), (v,Y) \in W \times \SS^f$ such that $(u,X) \preceq (v,Y)$. 
Set $u=u_0u_1$ and $v=v_0v_1$, where $u_0,v_0 \in W_T$, and $u_1,v_1$ are $(T,\emptyset)$-minimal.
Let $X_0 = T \cap u_1 W_X u_1^{-1}$ and $Y_0=T \cap v_1 W_Y v_1^{-1}$. 
Then $\pi_T(u,X) = (u_0,X_0)$ and $\pi_T(v,Y) = (v_0,Y_0)$. 
Let $w=v^{-1}u$, and let $w_0=v_0^{-1}u_0$. 
Since $(u,X) \preceq (v,Y)$, we have $X \subset Y$, $w \in W_Y$, and $w$ is $(\emptyset,X)$-minimal. 
We should show that $X_0 \subset Y_0$, $w_0 \in W_{Y_0}$, and $w_0$ is $(\emptyset,X_0)$-minimal. 
We argue by induction on the length of $w$. 
It is easily checked that, if $w=1$, then $u_0=v_0$ (thus $w_0=1$), $u_1=v_1$, and $X_0 \subset Y_0$, thus $\pi_T(u,X) \preceq \pi_T(v,Y)$. 
So, we may assume that $\lg(w) \ge 1$ plus the induction hypothesis.

\bigskip\noindent
We write $w=sw'$, where $s \in Y$, $w' \in W_Y$, and $\lg(w') = \lg(w)-1$. 
Let $v'=vs$. 
The element $(v')^{-1}u = w'$ lies in $W_Y$ and is $(\emptyset,X)$-minimal (since $w$ is), thus $(u,X) \preceq (v',Y)$. 
Set $v'=v_0' v_1'$, where $v_0' \in W_T$ and $v_1'$ is $(T, \emptyset)$-minimal, and set $Y_0' = T \cap (v_1') W_Y (v_1')^{-1}$.
By the induction hypothesis, we have $(u_0,X_0) = \pi_T(u,X) \preceq \pi_T(v',Y) = (v_0',Y_0')$. 
Set $w_0' = (v_0')^{-1}u_0$. 
Then $X_0 \subset Y_0'$, $w_0' \in W_{Y_0'}$, and $w_0'$ is $(\emptyset,X_0)$-minimal.

\bigskip\noindent
Suppose that $v_1s$ is $(T,\emptyset)$-minimal. 
Then $v_0'=v_0$ and $v_1'=v_1s$. 
Moreover, it is easily seen that, in that case, $Y_0 = Y_0'$ (thus $X_0 \subset Y_0$) and $w_0=w_0'$ (thus $w_0 \in W_{Y_0}$ and $w_0$ is $(\emptyset,X_0)$-minimal). 
Hence, $\pi_T(u,X) \preceq \pi_T(v,Y)$.

\bigskip\noindent
Suppose now that $v_1s$ is not $(T,\emptyset)$-minimal. 
We have $\lg(v_1s)>\lg(v_1)$, otherwise $v_1s$ would be $(T,\emptyset)$-minimal since $v_1$ is. 
Furthermore, by Proposition 2.9, there exists $t \in T$ such that $\lg(tv_1s) < \lg (v_1s)$. 
We also have $\lg(tv_1) > \lg(v_1)$, since $v_1$ is $(T,\emptyset)$-minimal. 
By Theorem~2.8, these inequalities imply that $tv_1=v_1s$. 
Then $v_0'=v_0t$, $v_1'=v_1$, thus $Y_0=Y_0'$ and $w_0=tw_0'$. 
A first consequence of this is that $X_0 \subset Y_0=Y_0'$ and $w_0 \in W_{Y_0}$ (since $w_0' \in W_{Y_0}$ and $t = v_1 s v_1^{-1} \in T \cap v_1 W_Y v_1^{-1} = Y_0$). 
It remains to prove that $w_0$ is $(\emptyset,X_0)$-minimal. 
Suppose not. 
Then we have $\lg(w_0) = \lg(tw_0') > \lg(w_0')$, otherwise $w_0$ would be $(\emptyset,X_0)$-minimal since $w_0'$ is. 
By Proposition~2.9, there exists $x \in X_0$ such that $\lg(t w_0' x)< \lg(tw_0')$. 
We also have $\lg( w_0' x) > \lg(w_0')$ since $w_0'$ is $(\emptyset,X_0)$-minimal. 
By Theorem~2.8, it follows that $t w_0' = w_0' x = w_0$. 
Hence, 
\[\begin{array}{lc}
&x=(w_0')^{-1} t (w_0') = u_0^{-1} (v_0') t (v_0')^{-1} u_0 \in W_{X_0} = W_T \cap u_1 W_X u_1^{-1}\\
\Rightarrow \quad & u^{-1} v_0 t v_0^{-1} u = u^{-1} v s v^{-1} u = w^{-1}sw \in W_X\\
\Rightarrow \quad & sw W_X = w' W_X = w W_X\,.
\end{array}\]
This contradicts the fact that $w$ is $(\emptyset,X)$-minimal (recall that $\lg(w') < \lg(w)$). 
So, $w_0$ is $(\emptyset,X_0)$-minimal.
We conclude that $\pi_T(u,X) \preceq \pi_T(v,Y)$.
\qed

\bigskip\noindent
Let $\bar \iota_T : \BSal(\Gamma_X) = \Sal(\Gamma_X)/W_X \to \BSal(\Gamma) = \Sal(\Gamma)/W$ denote the map induced by $\iota_T: \Sal(\Gamma_T) \to \Sal(\Gamma)$.

\bigskip\noindent
{\bf Lemma 5.2.}
{\it Let $T\subset S$.
Then $\bar \iota_T :\BSal(\Gamma_T) \to \BSal(\Gamma)$ is an embedding.}

\bigskip\noindent
{\bf Proof.}
Let $x,y$ be two points in $\Sal(\Gamma_T) \subset \Sal(\Gamma)$ that belong to the same $W$-orbit.
Let $w \in W$ such that $y = w\,x$.
The point $x$ (resp. $y$) lies in the interior of some cell $\B(u,X)$ (resp. $\B(v,Y)$) of $\Sal(\Gamma_T)$.
Since $w$ sends cells to cells, we shall have $w\, \B(u,X) = \B(w\,u,X) = \B(v,Y)$, thus $X=Y$ and $v=wu$. 
Hence, $w=vu^{-1} \in W_T$ (since $u,v \in W_T$).
\qed

\bigskip\noindent
Consider the cellular decomposition of $\BSal(\Gamma)$ described in Subsection 3.3.
Let $T$ be a subset of $S$.
Observe that, for all $X \in \SS^f_T$, the map $\bar\iota_T$ sends the cell $\bar \B (X)$ of $\BSal(\Gamma_T)$ homeomorphically to the cell $\bar \B (X)$ of $\BSal(\Gamma)$.
Observe also that, for all $s \in T$, the map $\bar\iota_T$ preserves the orientation of $\bar a_s = \bar \B(\{s\})$.
Hence,

\bigskip\noindent
{\bf Lemma 5.3.}
{\it Let $T$ be a subset of $S$.
Then the homomorphism $(\bar \iota_T)_* : \pi_1(\BSal(\Gamma_T), x_0) = A_{\Gamma_T} \to \pi_1(\BSal(\Gamma), x_0) = A_\Gamma$ coincides with the natural homomorphism $A_{\Gamma_T} \to A_\Gamma$ which sends $\sigma_s$ to $\sigma_s$ for all $s \in T$.}

\bigskip\noindent
{\bf Theorem 5.4}
(Van der Lek \cite{Lek1}).
{\it Let $T$ be a subset of $S$.
Then the natural homomorphism $A_{\Gamma_T} \to A_\Gamma$ which sends $\sigma_s$ to $\sigma_s$ for all $s \in T$ is injective. 
In other words, the pair $(A_T,\Sigma_T)$ is an Artin system of $\Gamma_T$.}

\bigskip\noindent
{\bf Proof.}
We have the following commutative diagram, where the lines are exact sequences. 
\[\begin{array}{ccccccccc}
1 & \to  & CA_{\Gamma_T} & \longrightarrow & A_{\Gamma_T} & \longrightarrow & W_{\Gamma_T} & \to & 1\\
  &      & \downarrow    &                 & \downarrow   &                 & \downarrow\\
1 & \to  & CA_{\Gamma}   & \longrightarrow & A_{\Gamma}   & \longrightarrow & W_{\Gamma}   & \to & 1
\end{array}\]
By Theorem 5.1, the homomorphism $CA_{\Gamma_T} \to CA_\Gamma$ has a retraction, thus it is injective. 
The homomorphism $W_{\Gamma_T} \to W_\Gamma$ is injective by \cite{Bourb1}.
We conclude by the five lemma that the homomorphism $A_{\Gamma_T} \to A_\Gamma$ is injective.
\qed

\bigskip\noindent
{\bf Theorem 5.5}
(Godelle, Paris \cite{GodPar1}).
{\it Let $T$ be a subset of $S$.
If $\Sal(\Gamma)$ is an Eilenberg MacLane space, then $\Sal(\Gamma_T)$ is also an Eilenberg MacLane space.}

\bigskip\noindent
{\bf Proof.}
By Theorem  5.1, for all $n \ge 1$, the homomorphism $(\iota_T)_*:\pi_n(\Sal(\Gamma_T),x(1)) \to \pi_n(\Sal(\Gamma),x(1))$ has a retraction $\pi_n(\Sal(\Gamma),x(1)) \to \pi_n(\Sal(\Gamma_T),x(1))$, hence it is injective.
Assume that $\Sal(\Gamma)$ is an Eilenberg MacLane space.
Then $\pi_n(\Sal(\Gamma),x(1))=\{1\}$ for all $n \ge 2$, thus, by the above, $\pi_n(\Sal(\Gamma_T),x(1)) = \{1\}$ for all $n \ge 2$, therefore $\Sal(\Gamma_T)$ is an Eilenberg MacLane space.
\qed

\bigskip\noindent
{\bf Theorem 5.6}
(Ellis, Sköldberg \cite{EllSko1}).
{\it Let $s,t \in S$ such that $m_{s,t} = \infty$.
Set $T=S \setminus \{s\}$ and $R=S\setminus \{t\}$.
If $\Sal(\Gamma_T)$ and $\Sal(\Gamma_R)$ are both Eilenberg MacLane spaces, then $\Sal(\Gamma)$ is an Eilenberg MacLane space.}

\bigskip\noindent
{\bf Proof.}
Let $X \in \SS^f$.
If $s \in X$, then $t \not\in X$, thus $X \in \SS^f_R$.
Similarly, if $t \in X$, then $s \not\in X$, thus $X \in \SS^f_T$.
Note that, if $s \not\in X$ and $t \not\in X$, then $X \in \SS_{T \cap R}^f$.
By Lemma 5.2, $\BSal( \Gamma_T)$ and $\BSal(\Gamma_R)$ are subcomplexes of $\BSal(\Gamma)$.
By the above, we have $\BSal(\Gamma_T) \cup \BSal(\Gamma_R) = \BSal(\Gamma)$ and $\BSal(\Gamma_T) \cap \BSal (\Gamma_R) = \BSal(\Gamma_{T \cap R})$.
By Lemma 5.3, the homomorphisms $\pi_1(\BSal (\Gamma_{T \cap R}), x_0) \to \pi_1(\BSal (\Gamma_{T}), x_0)$ and $\pi_1(\BSal (\Gamma_{T \cap R}), x_0) \to \pi_1(\BSal (\Gamma_{R}), x_0)$ are injective.
The complexe $\BSal(\Gamma_T)$ is an Eilenberg Maclane space, since $\Sal(\Gamma_T)$ is a covering of $\BSal(\Gamma_T)$ and, by hypothesis, $\Sal(\Gamma_T)$ is an Eilenberg MacLane space.
Similarly, $\BSal(\Gamma_R)$ is an Eilenberg MacLane space.
Furthermore, by applying Theorem 5.5 to $\Gamma_{T \cap R}$ and $\Gamma_T$ we get that $\Sal(\Gamma_{T \cap R})$ is an Eilenberg MacLane space, thus $\BSal(\Gamma_{T \cap R})$ is an Eilenberg MacLane space.
By Theorem~2.4, it follows that $\BSal(\Gamma)$ is an Eilenberg MacLane space, thus $\Sal(\Gamma)$ is an Eilenberg MacLane space.
\qed

\bigskip\noindent
Recall that a subset $T$ of $S$ is said to be {\it free of infinity} if $m_{s,t} \neq \infty$ for all $s,t \in T$.
We denote by $\SS^{<\infty}$ the set of subsets of $S$ that are free of infinity.
Note that $\SS^f \subset \SS^{< \infty}$.
Recall also that $\Gamma$ is said to be of {\it FC type} if $\SS^f = \SS^{< \infty}$.

\bigskip\noindent
{\bf Corollary 5.7.}
(Ellis, Sköldberg \cite{EllSko1}).
{\it If $\Sal(\Gamma_T)$ is an Eilenberg MacLane space for all $T \in \SS^{<\infty}$, then $\Sal(\Gamma)$ is an Eilenberg MacLane space.}

\bigskip\noindent
{\bf Proof.}
We argue by induction on $|S|$.
If $|S|=1$, then $S$ is free of infinity, thus $\Sal(\Gamma)=\Sal(\Gamma_S)$ is an Eilenberg MacLane space.
More generally, if $S$ itself is free of infinity, then $\Sal(\Gamma)=\Sal(\Gamma_S)$ is an Eilenberg MacLane space.
Assume that $|S| \ge 2$ and $S$ is not free of infinity, plus the induction hypothesis.
Let $s,t \in S$ such that $m_{s,t} = \infty$.
Set $T=S \setminus \{s\}$ and $R = S \setminus \{t\}$.
By the induction hypothesis, $\Sal(\Gamma_T)$ and $\Sal(\Gamma_R)$ are Eilenberg MacLane spaces.
By Theorem 5.6, it follows that $\Sal(\Gamma)$ is an Eilenberg MacLane space.
\qed

\bigskip\noindent
{\bf Corollary 5.8.}
(Charney, Davis \cite{ChaDav1}).
{\it The complex $\Sal(\Gamma)$ is an Eilenberg MacLane space if $\Gamma$ is a Coxeter graph of FC type.}

\bigskip\noindent
{\bf Proof.}
By Theorem 4.10, $\Sal(\Gamma_T)$ is an Eilenberg MacLane space for all $T \in \SS^f$.
By definition, $\SS^f = \SS^{<\infty}$, thus, by Corollary 5.7, $\Sal(\Gamma)$ is an Eilenberg MacLane space.
\qed

%%%%%%%%%%%%%%%%%%%%%%%
%%%%%BIBLIOGRAPHIE%%%%%
%%%%%%%%%%%%%%%%%%%%%%%

%%%%%%%%%%

\bigskip\bigskip\noindent
{\bf Luis Paris,}

\smallskip\noindent 
Université de Bourgogne, Institut de Mathématiques de Bourgogne, UMR 5584 du CNRS, B.P. 47870, 21078 Dijon cedex, France.

\smallskip\noindent
E-mail: {\tt lparis@u-bourgogne.fr}

%%%%%%%%%%

\end{document}